%5.02.97
%length 85722
%30 pages
\input amstex
\documentstyle{amsppt}
\magnification=1200
\hsize=16.5truecm
\vsize=23.3truecm
%\voffset=1.5truecm

\catcode`\@=11
\redefine\logo@{}
\catcode`\@=13
\define\gt{\Gamma_t}
\define\veps{\varepsilon}
\define\vth{\vartheta}

\define\al{\operatorname{Lift}}
\define\dv{\operatorname{Div}}

\define \bn{\Bbb N}
\define \bz{\Bbb Z}
\define \bq{\Bbb Q}
\define \br{\Bbb R}
\define \bc{\Bbb C}
\define \bh{\Bbb H}

\define\gi{\Gamma_{\infty}}

\define\pd#1#2{\dfrac{\partial#1}{\partial#2}}

\input epsf.tex

\TagsOnRight
%\NoBlackBoxes
\document

\topmatter
\title
Commutator coverings of Siegel threefolds
\endtitle

\author
V.  Gritsenko \footnote{Supported by RIMS of Kyoto University
\hfill\hfill}
and K. Hulek$^1$
\endauthor
%\date 5 February,\ 1997
%,${}\quad{}$ alg-geom/???????
%\enddate

\address
St. Petersburg Department of Steklov Mathematical Institute,
\newline
${}\hskip 8pt $ Fontanka 27, 191011 St. Petersburg, Russia
\endaddress
\email
gritsenk\@kurims.kyoto-u.ac.jp (till 28.02.97);
gritsenk\@mpim-bonn.mpg.de
\endemail

\address
Institut f\"ur Mathematik, Universit\"at Hannover, Postfach 6009,
\newline
${}\hskip 9pt $ D-30060 Hannover, Germany
\endaddress

\email
hulek\@math.uni-hannover.de
\endemail

 \abstract
We investigate the existence and non-existence of modular forms of low weight
with a character with respect to the paramodular group $\Gamma_t$ and discuss
the resulting geometric consequences. Using an advanced version of Maa\ss\
lifting one can construct many examples of such modular forms and in particular
examples of weight 3 cusp forms. Consequently we find many abelian coverings of
low degree of the moduli space ${\Cal A}_t$ of $(1,t)$-polarized abelian
surfaces which are not unirational. We also determine the commutator subgroups
of the paramodular group $\Gamma_t$ and its degree 2 extension $\Gamma^+_t$.
This has applications for the Picard group of the moduli stack
${\underline{\Cal A}}_t$.
Finally we prove non-existence theorems for low weight modular forms.
As one of our main results we obtain the theorem that the
maximal abelian cover ${\Cal A}_t^{com}$ of ${\Cal A}_t$ has geometric genus $0$
if and only if $t=1$, $2$, $4$ or $5$.
We also prove that ${\Cal A}_t^{com}$ has geometric genus $1$ for
$t=3$ and $7$.
\endabstract

\rightheadtext
{Commutator coverings of Siegel threefolds}
\leftheadtext{V. Gritsenko and  K. Hulek}
\endtopmatter

\document

\document

\head
Introduction
\endhead
The main theme of this paper are cusp forms of small weight with
respect to the paramodular group $\Gamma_t$ and Siegel modular
threefolds. Since cusp forms of weight 3 define canonical differential
forms on the moduli space $\Cal A_t$ of $(1,t)$--polarized abelian
surfaces, the existence or non-existence of such forms has important
geometric consequences.

For $t\ge 1$ the paramodular group $\Gamma_t$ is not a maximal discrete
group. One can add to it a number of exterior involutions to obtain a
maximal normal extension $\Gamma_t^*$ such that $\Gamma_t^*/\Gamma_t$
is a product of $\bz_2$ components. We then have the following tower
of Siegel threefolds
$$
\CD
\Cal A_t^{com}@>>>\Cal A^{(l)}@>>>
\Cal A_t@>>> \Cal A^{(r)}@>>>\Cal A_t^{*}\\
@|@|@|@|@|\\
\Gamma_t'\setminus \bh_2@>>>\Gamma^{(l)}\setminus \bh_2@>>>
\Gamma_t\setminus \bh_2
@>>>\Gamma^{(r)}\setminus \bh_2@>>>
\Gamma_t^{*}\setminus \bh_2
\endCD
\tag{0.1}
$$
where $\Gamma_t'\subset \Gamma^{(l)}\subset \Gamma_t
\subset \Gamma^{(r)}\subset \Gamma_t^*$
and $\Gamma_t'$ is the commutator subgroup of $\Gamma_t$.
We call  these varieties
{\it commutative neigbours} of ${\Cal A}_t$.
They have very interesting properties.
Neighbours $\Cal A^{(r)}$ of the spaces
$\Cal A_t$ to the right, i.e. finite quotients, were studied in \cite{GH2}.
The coverings in \thetag{0.1} to the left of ${\Cal A}_t$ are  also
galois with a finite abelian Galois group.
This paper is devoted to neighbours $\Cal A^{(l)}$ to the left
of ${\Cal A}_t$.

Recall that the geometric genus of $\Cal A_t$ can be zero only for $20$
exceptional polarizations 
$t=1, \dots,\ 12$, $14$, $15$, $16$, $18$, $20$, $24$, $30$, $36$
(see \cite{G1}).
Hypothetically, for all of them $\Cal A_t$ is rational or unirational
(see \cite{GP} for an announcement of some  results in this
direction  and \cite{G4}). Hence no weight $3$ cusp forms should exist
for these values of $t$ 
(for an easy proof of this fact for $t\le 8$ see Corollary 3.3).
On the other hand we shall construct
many examples of cusp forms of small weight $(k\le 3)$ with respect to
$\Gamma_t$ with a character. It follows that $\Cal A_t$ usually has a
double modular covering with positive geometric genus. One of the
main results of this paper is

\proclaim{Theorem 0.1}
Let ${\widetilde{\Cal A}}_t^{com}$ be a smooth projective model of
the maximal abelian covering $\Cal A_t^{{com}}=\Gamma_t'\setminus
\bh_2$  of ${\Cal A}_t$. Then
\roster
\item
$\ \widetilde{\Cal A}_t^{com}\ $ has geometric genus $0$ if and
only if $t=1$, $2$, $4$ or $5$;
\smallskip
\item
$\ {\widetilde{\Cal A}}_3^{com}\ $ and
$\ {\widetilde{\Cal A}}_7^{com}\ $ have geometric genus 1.
\endroster
\endproclaim

Therefore ${\widetilde{\Cal A}}_3^{{com}}$ and ${\widetilde{\Cal
A}}_7^{{com}}$ are candidates for Calabi--Yau varieties. Some other
examples of Siegel threefolds with genus $1$ can also be found in section 3
(see Theorem 3.1 and Corollary 3.5).

The paper is organized as follows.
In section 1 we collect examples of $\Gamma_t$-cusp forms with a
character of weight one, two  or three. The main ingredient for the
construction of these forms is an advanced version of the Maa\ss\
lifting of Jacobi forms to modular forms. This lifting result was
proved in \cite {GN} and goes back to \cite {G1}. As a consequence of
this construction the modular forms which we obtain are in fact forms
with respect to a bigger group $\Gamma_t^+$ which is an extension of
$\Gamma_t$ of order 2.  Characters of $\Gamma_t$ define coverings of
$\Cal A_t$ (see definition 1.5). As a consequence we obtain many
covering spaces of $\Cal A_t$ for which we can prove that they are not
unirational, resp. in some cases we can also show that they have
Kodaira dimension $\ge 1$.  This result is contained in Corollary 1.6.
Note that these spaces are covering spaces of rational or unirational
moduli spaces $\Cal A_t$ of degree $2$, $3$ or $4$.

In section 2 we determine, partly using the existence of characters
from section 1, the commutator subgroups of $\Gamma_t$ and
$\Gamma_t^+$ (Theorem 2.1). This result could be classical (in the
case of $Sp_4(\Bbb Z)$ see \cite {R} and \cite {Ma}), but to our
knowledge it does not exist in the literature. We also note some
consequences of this result for the Picard group of the moduli stack
$\underline{\Cal A}_t$ of $(1,t)$--polarized abelian surfaces. In
particular the quotient of $\Gamma_t$ by its commutator subgroup gives
the torsion part of the Picard group (see Proposition 2.3).

In section 3 we study the maximal abelian covering $\Cal A_t^{{com}}$
of $\Cal A_t$. The crucial technical ingredient in section 3 is an
estimate on the vanishing order of Jacobi cusp forms (Proposition
3.2). This result seems interesting in its own right. As an easy
corollary one obtains that there exist no weight 3 cusp forms for
$\Gamma_t$ for $t\le 8$ (Corollary 3.3).  Among the exceptional
polarizations the  polarization $(1,3)$ looks very
interesting.
There are twelve commutative neighbours of the threefold
$\Cal A_3$ to the left.
Six of them have geometric genus $1$ and
the others have genus $0$.
We also found an example of a Siegel modular threefold,
namely a 3:1-covering of $\Cal A_6$, which has geometric genus
$2$ (see Theorem 3.1).
\medskip

\noindent {\it Acknowledgment. } It is a pleasure for us to thank RIMS of Kyoto
University for support. We were very impressed by the hospitality of
RIMS and its stimulating working atmosphere. We would also like to thank
the DFG for partial support under 436 RUS 17/108/95.
The second author would like to thank W. Fulton for discussions concerning
\cite{Mu1} during a visit in Cicago.

\medskip

\head
\S 1. Roots of order $2$, $3$, $4$, $6$ and $12$ of
Siegel modular forms
\endhead

In this section we construct roots of Siegel modular forms with
respect to the paramodular group $\gt$ $(t\in \bn)$.  They will be
again modular forms with respect to the full paramodular group (no
level structure involved) but with a character $\chi: \gt \to \bc^*$.

We will use two realizations of the paramodular group.  Let us fix a
free $\bz$-module $L=\Bbb Z e_1\oplus \Bbb Z e_2 \oplus \Bbb Z e_3
\oplus \Bbb Z e_4$.  We define a skew-symmetric form $W_t$ on $L$ as
follows
$$
W_t(x,y)e_1\wedge e_2\wedge e_3\wedge e_4= -x\wedge y\wedge w_t,\qquad
w_t=te_1\wedge e_3+ e_2\wedge e_4.
$$
The group
$$
\widetilde{\Gamma}_t =Sp(W_t, \bz)= \{g:L\to L\ |\  W_t(gx,gy)=W_t(x,y)\}
$$
is called {\it the integral paramodular group of level} $t$.  It is
easy to check that
$$
\text{if }\ g=(g_{ij})\in \widetilde{\Gamma}_t,
\quad\text{then }\ g_{12}\equiv g_{14}\equiv g_{32}\equiv g_{34}
\equiv 0
\hbox{ mod }t.
$$
The paramodular group $\widetilde{\Gamma}_t$ is conjugate to a
subgroup of the rational symplectic group $Sp_4(\bq)$:
$$
\Gamma_t:=I_t \widetilde{\Gamma}_t I_t^{-1}=
\left\{\pmatrix
*  & t* & * & *\\
*  & * & * & t^{-1}*\\
*  & t * & * & *\\
t* & t*& t* & *
\endpmatrix
 \in Sp_4 (\bq)\,|\, \text{ all $*$ are integral}\, \right\}
$$
where $I_t=\hbox{diag}\,(1,1,1,t)$.  This is the second realization of
the paramodular group.  The quotient
$$
\Cal A_t=\Gamma_{t}\setminus \bh_2,
$$
where $\bh_2$ is the Siegel upper half-plane of genus 2, is isomorphic
to the moduli space of abelian surfaces with a polarization of type
$(1,t)$.

A modular form of weight $k$ with respect to $\Gamma_t$ with a
character $\chi: \Gamma_t\to \bc^*$ is a holomorphic function on
$\bh_2$ which for arbitrary $M\in \gt$ satisfies the functional
equation
$$
F(M<Z>)=\chi(M)\hbox{det}(CZ+D)^{k}F(Z) \qquad M=\left(\smallmatrix
A&B\\C&D\endsmallmatrix\right) \in \gt.
$$
We denote the space of such modular (cusp) forms by $\frak M_{k}(\gt,
\chi)$ ($\frak N_{k}(\gt, \chi)$ respectively).  Here we admit a
character of the paramodular group $\Gamma_t$ in order to construct
roots of certain orders of modular forms with respect to $\Gamma_t$
with trivial character. We remark that no level structure is involved
in our considerations.

The group $\gt$ is not a maximal discrete group acting on $\bh_2$ if
$t\ne 1$. It has a normal extension $\Gamma_t^*$ such that
$\Gamma_t^*/\gt\cong (\bz/2)^{\nu(t)}$ where $\nu(t)$ is the number
of distinct prime divisors of $t$.  The modular forms which we
construct in this paper are forms with respect to a double extension
of $\gt$:
$$
\gt^+=<\gt, V_t> \qquad\text{where }\quad
V_t={\tsize\frac 1{\sqrt{t}}}\left(\smallmatrix 0 & t & 0 & 0\\ 1 & 0
& 0 & 0\\ 0 & 0 & 0 & 1\\ 0 & 0 & t & 0
\endsmallmatrix\right).
\tag{1.1}
$$
The double quotient
$$
\Cal A_t\overset{2:1}\to\longrightarrow
\Cal A_t^+=\gt^+\setminus \bh_2
$$
of $\Cal A_t$ can be interpreted as a moduli space of
lattice-polarized $K3$ surfaces for arbitrary $t$ or as the moduli
space of Kummer surfaces of $(1,p)$-polarized abelian surfaces for a
prime $t=p$ (see \cite{GH2}).

The main tool of the construction of modular forms of type
$\root{n}\of{F(Z)}$, where $F(Z)\in \frak M_{k}(\gt)$ is a new variant
of the lifting of Jacobi forms to Siegel modular forms proposed in
\cite{GN, Theorem 1.12}.
The datum for this lifting is a Jacobi form $\phi_{k,R}(\tau,z)$ of
integral weight $k$ and index $R$ ($R\in \bn/2$) together with a
character of the full Jacobi group. We define Jacobi forms as modular
forms with respect to the maximal parabolic subgroup $\Gamma_\infty$
of $Sp_4(\bz)=\Gamma_1$ which fixes the line spanned by $e_2$, i.e.
$$
\gi=
\left\{\left(\smallmatrix *&0&*&*\\
                *&*&*&*\\ *&0&*&*\\ 0&0&0&*
\endsmallmatrix\right)
\in Sp_4(\Bbb Z)\right\}.
$$
Let $k$ and $R$ be integral or half-integral.  We choose one of the
holomorphic square roots $\sqrt {\hbox{det} (Z)}$ by the condition
$\sqrt{\hbox{det}(Z/i)}>0$ for $Z=i{\bold 1}$.

We call a holomorphic function $\phi(\tau,z)$ on $\bh_1\times \bc$
{\it a Jacobi form of weight $k$ and index $R$ with a multiplier
system (or a character)} $v:\Gamma_\infty\to \bc^{*}$ if the function
$$
\widetilde{\phi}(Z):=\phi(\tau ,z)\exp{(2\pi i R\omega )},
\qquad
Z=\pmatrix \tau&z\\z&\omega\endpmatrix\in \bh_2,
$$
satisfies the functional equation
$$
\widetilde{\phi} (M<Z>)=v(M)\hbox{det}(CZ+D)^{k}\widetilde{\phi}(Z)
\qquad\qquad \text{for any }
M=\left(\smallmatrix A&B\\C&D\endsmallmatrix\right)\in
\Gamma_{\infty}
$$
and $\phi$ has a Fourier expansion of type
$$
\phi(\tau ,z)=\sum\Sb n,\,l\\
\vspace{0.5\jot} 4Rn-l^2\ge 0\endSb
f(n,l) \exp{(2\pi i(n\tau +lz))}
$$
where the summation is taken over $n$ and $l$ from some free
$\bz$-modules depending on $v$.  The condition $f(n,l)=0$ unless
$4Rn-l^2\ge 0$ is equivalent to the holomorphicity of $\phi$ at
infinity.  The form $\phi(\tau ,z)$ is called a Jacobi cusp form if
$f(n,l)=0$ unless $4Rn-l^2>0$.  We call the number $ 4Rn-l^2$ {\it the
norm} of the index of the Fourier coefficient $f(n,l)$.  We denote the
finite dimensional space of such Jacobi forms (cusp forms) by
$J_{k,R}(v)$ ($J_{k,R}^{c}(v)$ respectively).

If the function $\phi$ has a Fourier expansion of type
$$
\phi(\tau ,z)=\sum\Sb n \ge 0,\,l\endSb
f(n,l) \exp{(2\pi i(n\tau +lz))},
\tag{1.2}
$$
then we call it a {\it weak Jacobi form}. The space $J_{k,R}^{w}(v)$
of all such forms is again finite dimensional.

We remark that the group $\gi/\{\pm {\bold 1}_2\}\cong
SL_2(\bz)\ltimes H(\bz)$, where $H(\bz)$ is the integer Heisenberg
group, is usually called the Jacobi group.  We use the following
embeddings of $SL_2(\bz)$ and $H(\bz)$ into $\Gamma_{\infty}$
$$
\hskip-1pt i_\infty(
\left(\smallmatrix a&b\\c&d\endsmallmatrix\right)):=
\left(\smallmatrix
a&0&b&0\\ 0&1&0&0\\ c&0&d&0\\ 0&0&0&1
\endsmallmatrix\right),\qquad
H(\Bbb Z)\cong\left\{\lambda,\mu,\kappa\in \bz\ |\ [\lambda,\mu;
 \kappa]=
\left(\smallmatrix
1&0&0&\hphantom{-}\mu\\
\lambda&1&\mu&\hphantom{-}\kappa\\
0&0&1&{-}\lambda\\ 0&0&0&\hphantom{-}1
\endsmallmatrix\right)
\right\}.
\tag{1.3}
$$
Here we admit Jacobi forms of {\it half-integral} indices.  This is
the only difference between the definition of the Jacobi forms given
above and the definition of \cite{EZ}.

The best known examples of forms with half-integral indices are the
Jacobi triple product and the quintiple product.  The Jacobi
theta-series, which was, maybe, the first example of a Jacobi modular
form, is defined as
$$
\vartheta(\tau ,z)=\hskip-2pt\sum\Sb n\equiv 1\, mod\, 2 \endSb
\,(-1)^{\frac{n-1}2}
\exp{(\frac{\pi i\, n^2}{4} \tau +\pi i\,n z)}=
\sum_{m\in \bz}\,\biggl(\frac{-4}{m}\biggr)\, q^{{m^2}/8}\,r^{{m}/2},
\tag{1.4}
$$
where $q=\exp{(2\pi i \tau)}$, $r=\exp{(2\pi i z)}$ and

$$
\biggl(\frac{-4}{m}\biggr)=\cases \pm 1 &\text{if }
m\equiv \pm 1\ \hbox{mod}\ 4\\
\hphantom{\pm}0 &\text{if }
m\equiv \hphantom{\pm} 0 \ \hbox{mod}\ 2.
\endcases
$$
This is a Jacobi form of weight $1/2$ and index $1/2$, i.e.  an
element of $J_{\frac 1{2}, \frac1{2}}(v_\eta^3\times v_H)$.  The
multiplier system $v_\eta^3\times v_H$ is induced by the
$SL_2(\bz)$-multiplier system $v_\eta$ of the Dedekind $\eta$-function
and the following character of the integral Heisenberg group
$$
v_H([\lambda,\mu;\kappa]):=(-1)^{\lambda+\mu+\lambda\mu+\kappa}.
\tag{1.5}
$$
We recall the famous Jacobi triple product formula
$$
\vartheta(\tau ,\,z)=
-q^{1/8}r^{-1/2}\prod_{n\ge 1}\,(1-q^{n-1} r)(1-q^n r^{-1})(1-q^n).
\tag{1.6}
$$

The quintiple product is a Jacobi form of weight $1/2$ and index $3/2$
$$
\vartheta_{3/2}(\tau ,z)=\sum_{n\in \bz}
\biggl(\frac{12}n\biggl)\,q^{{n^2}/{24}}r^{{n}/2}
\in J_{\frac 1{2}, \frac 3{2}}(v_\eta\times v_H)
\tag{1.7}
$$
where
$$
\biggl(\frac{12}n\biggl)=
\cases
\hphantom{-}1\  \text{ if }\  n\equiv \pm 1 \operatorname{mod} 12\\
-1\ \text{ if }\ n\equiv \pm 5 \operatorname{mod} 12\\
\hphantom{-}0\  \text{ if }\  (n,12)\ne 1.
\endcases
$$
It is related to the Jacobi theta-series by the formula
$$
\vartheta_{3/2}(\tau ,z)
=\frac{\eta(\tau )\vartheta(\tau , 2z)} {\vartheta(\tau , z)}.
$$

Many examples of Jacobi forms of half-integral indices were
constructed in \cite{GN}.  One can, for instance, define Hecke
operators of different types on Jacobi forms.  The first type of Hecke
operator is the map
$$
\phi(\tau, z)\mapsto \phi(\tau, az) \qquad(a\in \bn).
$$
For the Jacobi theta-series we have
$$
\vth_a:= \vartheta(\tau,az)
\in J_{\frac1{2}, \frac 1{2}{a^2}}(v_\eta^3\times v_H^a).
$$

The second type of Hecke operator which we shall use here are the
operators $T^{(Q)}_-(m)$ acting on Jacobi forms in
$J_{k,R}(v_\eta^{24/Q}\times v_H^\varepsilon)$, where $Q=1$, $2$, $3$,
$4$, $6$ or $12$, $\varepsilon\equiv 2R\ \hbox{mod}\ 2$, the weight
$k$ is an integer and $m\equiv 1\ \hbox{mod}\ Q$.  We assume for
simplicity that $QR\in \bn$ as well. By definition
$$
\widetilde{\phi}|_k T^{(Q)}_-(m)(Z)=
m^{2k-3}\sum\Sb ad=m\\ \vspace {0.5\jot} b\, mod\, d\endSb
d^{-k}v_\eta^{24/Q}(\sigma_a)\phi(\frac{a\tau +bQ}d,\,az)
\exp{(2\pi i mR\omega)}
\tag{1.8}
$$
where $\sigma_a\in SL_2(\bz)$ such that $\sigma_a\equiv
\left(\smallmatrix a^{-1}&0\\0&a\endsmallmatrix\right)
\hbox{ mod}\, Q$.
One can show that
$$
\phi|_k T^{(Q)}_-(m)\in
J_{k,mR}(v_\eta^{24/Q}\times v_H^\varepsilon)
$$
(see \cite{GN, Lemma 1.7}).  The lifting we shall use in this paper is
a new variant of the lifting construction of \cite{G1} proposed in
\cite{GN} (see also \cite{G3} for the case of orthogonal
groups).  All these liftings are advanced versions of the Maa\ss\
lifting.  Let us consider a Jacobi cusp form $
\phi\in J_{k,R}^{c}(v_\eta^{24/Q}\times v_H^{\varepsilon}),
$ where $Q$, $R$, and $\varepsilon$ are as above.  Then the function
$$
\al(\phi)(Z)
=\sum\Sb m\equiv 1 \, mod \,Q\\
\vspace{0.5\jot} m>0\endSb
m^{2-k}\bigl(\widetilde{\phi}|_k T_-^{(Q)}(m)\bigr)(Z)
\in \frak N_k(\Gamma_{t}, \chi_Q)
\tag{1.9}
$$
is a non-zero cusp form with respect to the paramodular group
$\Gamma_t$ with $t=QR\in \bn$, where
$$
\chi_Q:\Gamma_{t}\to \{\root{Q}\of{1}\}
$$
is a character of order $Q$ of $\Gamma_{t}$ induced by
$v_{\eta}^{24/Q}\times v_H^{\varepsilon}$ and $\veps\equiv \frac{2t}Q\
\hbox{mod}\ 2$.  I.e. we have
$$
\chi_Q|_{i_{\infty}(SL_2(\bz))}=v_\eta^{24/Q},\qquad
\chi_Q|_{H(\bz)}=v_H^{{2t}/Q}.
\tag{1.10a}
$$
Here $SL_2(\bz)$ and $H(\bz)$ are realized as subgroups of $\gi\subset
\gt$ (see \thetag{1.3}) and
$$
\chi_Q([0,0;\frac{\kappa}{t}])=
\exp{(2\pi i \frac{\kappa}Q)}
\tag{1.10b}
$$
where $\{[0,0;{\kappa}/{t}] \,|\, \kappa\in \bz\}$ is the center of
the maximal parabolic subgroup $\Gamma_{t,\infty}=\gi(\bq)\cap
\Gamma_{t}$ of the group $\Gamma_{t}$ (\cite{GN, Theorem 1.12}).
Moreover $\al(\phi)$ satisfies the equation
$$
\al(\phi)(V_{t}<Z>)=\al(\phi)(Z)
$$
for the involution \thetag{1.1}. Thus $\al(\phi)$ is in fact a cusp
form with respect to the double extension $\Gamma_{t}^+$ of
$\Gamma_{t}$, i.e.
$$
\al(\phi)\in
\frak N_k(\Gamma_{t}^+, \chi_{Q,{(-1)}^k}),
\tag{1.11}
$$
Here $\chi_{Q,\pm}$ denotes a character of $\Gamma_{t}^+$ such that
$$
\chi_{Q,\pm}|\Gamma_t=\chi_Q, \qquad
\chi_{Q,\pm}(V_t)=\pm 1.
$$
We shall see in \S 2 that the last relation and
\thetag{1.10a--b} determine the character $\chi_{Q,\pm}$ uniquely.
In the sequel we shall sometimes add an additional index
$\chi_{Q,\pm}^{(t)}$ to indicate for which paramodular group the
character is defined.

We can consider the modular form $\al(\phi_{k,R})$ as a root of order
$Q$ of a modular form with respect to $\Gamma_{t}$ ($t=QR$) with a
trivial character.  This construction is specially interesting in the
case of small weights. For example, there are no modular forms of
weight $1$ with respect to $Sp_4(\bz)$ with trivial character.  It
follows from a well known result of Skoruppa, that there are also no
Jacobi forms of weight one with respect to the full Jacobi group (see
\cite{Sk}).  On the other hand we shall construct infinitely many
$\Gamma_t$-forms, which are indeed cusp forms, of weight $1$ with
characters of order $4$, $6$ and $12$. We will also construct several
series of $\Gamma_t$-cusp forms of weight $2$ and $3$ for all $t\ne
1,\,4,\,5$.  We remark that the first cusp form of weight $3$ with
respect to $\Gamma_t$ with trivial character exists for $t=13$ (see
\cite{G1}).

In \S 2 we shall determine the commutator of the groups $\Gamma_t$ and
$\Gamma_{t}^+$. For this we need the existence of certain characters
of $\Gamma_{t}^+$ which we can obtain from the following lemma.  The
non-trivial binary character of the modular group $\Gamma_1=Sp_4(\bz)$
was constructed in \cite{R} and \cite{Ma}, where it was proved that
the commutator subgroup of $Sp_4(\bz)$ has index $2$.

\proclaim{Lemma 1.1}There exist characters
$\chi_{Q,\pm}^{(t)}$ of order $Q$ of the group $\Gamma_t^+$ in the
following cases
$$
\chi^{(t)}_{2,\pm}\ \ \text{ \rm for  arbitrary }\ t,
\qquad
\chi^{(t)}_{3,\pm}\ \ \text{ \rm for }\ t\equiv 0\hbox{ \rm mod }3,
\qquad
\chi^{(t)}_{4,\pm}\ \ \text{ \rm for }\ t\equiv 0\hbox{ \rm mod }2.
$$
These characters are induced by the character $v_\eta^{{24}/Q}\times
v_H^{{2t}/Q}$.
\endproclaim
\demo{Proof} Since $\Gamma_t$ is a normal subgroup of $\Gamma_{t}^+$
of index $2$ there exists a character of $\Gamma_{t}^+$ which is
trivial on $\Gamma_t$ and has value $-1$ on $V_t$. Hence it is enough
to construct one character for $\Gamma_{t}^+$ of type $\chi_{Q,\pm}$.

We consider the Jacobi-Eisenstein series
$$
e_{4,1}(\tau,z)=1+(r^2+56r+126+56r^{-1}+r^{-2})q+\dots\in J_{4,1}
$$
of weight $4$ and index $1$ (see \cite{EZ, \S 2}). We set
$$
e_{4,m}=e_{4,1}|_4 T_-(m)\in J_{4,m}
$$
where $T_-(m)=T^{(1)}_-(m)$ is the Hecke operator \thetag{1.8}.  Using
the lifting construction \thetag{1.9} we obtain
$$
\al(e_{4,m}\eta^6\vth^2)\in \frak N_8(\Gamma_{2m+2}^+, \chi_{2,+}),
\qquad
\al(e_{4,m}\eta^3\vth^3)\in \frak N_7(\Gamma_{2m+3}^+, \chi_{2,-})
$$
where we put $e_{4,0}=1$ for $m=0$.  For $t=1$ the non-trivial
character of $\Gamma_1$ is the character of the cusp form
$\al(\eta^9\vth)$.  To obtain further characters we consider
$$
\align
\al(e_{4,m}\eta^3\vth_2)&\in \frak N_6(\Gamma_{4m+8}^+, \chi_{4,+}),\\
\al(e_{4,m}\eta^3\vth)&\in \frak N_6(\Gamma_{4m+2}^+, \chi_{4,+}),\\
\al(e_{4,m}\eta^2\vth^2)&\in \frak N_6(\Gamma_{3m+3}^+, \chi_{3,+}).
\endalign
$$
To complete the list of characters we add the character
$\chi_{4,+}$ of the group $\Gamma_4$.  This is the character of the
lifting $\al(\phi_{10,m}\vth^2)$, where
$\phi_{10,m}=\phi_{10,1}|_{10}T_-(m)\in J_{10,m}$.  The Jacobi form
$\phi_{10,1}=\eta^{18}\vth^2$ is the first Jacobi cusp form of index
$1$.

\hfill\hfill\qed
\enddemo
Our next aim is to construct $\Gamma_t$-cusp forms of small weights.
We shall explain later in this section why this is of geometric
interest.  To construct such forms we first define Jacobi cusp forms
of weights between $1$ and $2$ using the Jacobi theta-series.  The
next lemma is an extended version of Lemma 1.18 of \cite{GN}.

\proclaim{Lemma 1.2} Let $a,\, b,\,c,\,d\in \bn$.
Then
$$
\aligned
\vartheta_a\vartheta_b
&\in J_{1, \frac 1{2} ({a^2+b^2})}^{c}(v_\eta^6\times v_H^{a+b})\\
(\vartheta_{3/2})_a(\vartheta_{3/2})_b
&\in J_{1,\frac{3}{2}(a^2+b^2)}^{c}(v_\eta^2\times v_H^{a+b})\\
\vartheta_a\vartheta_b\vth_c
&\in J_{\frac3{2}, \frac 1{2} ({a^2+b^2+c^2})}^{c} (v_\eta^9\times
v_H^{a+b+c})\\ (\vth_{3/2})_a(\vth_{3/2})_b(\vth_{3/2})_c
&\in J_{\frac3{2}, \frac 3{2} ({a^2+b^2+c^2})}^{c}
(v_\eta^3\times v_H^{a+b+c})\\
\vartheta_a\vartheta_b\vth_c\vth_d
&\in J_{2, \frac 1{2} ({a^2+b^2+c^2+d^2})}^{c}
(v_\eta^{12}\times v_H^{a+b+c+d})\\
(\vth_{3/2})_a(\vth_{3/2})_b(\vth_{3/2})_c(\vth_{3/2})_d &\in J_{2,
\frac 3{2} ({a^2+b^2+c^2+d^2})}^{c} (v_\eta^{4}\times v_H^{a+b+c+d})
\endaligned
\
\aligned
{}&\text{if }\ \tsize\frac{ab}{(a,b)^2}\ \text{ is even}\\
{}&\text{if}\ \left(\tsize\frac{ab}{(a,b)^2}, 6\right)\ne 1\\
{}&\text{if }\ \tsize\frac{abc}{(a,b,c)^3}\ \text{ is even}\\
{}&\text{if }\ \left(\tsize\frac{abc}{(a,b,c)^3},6\right)\ne 1\\
{}&\text{if }\ \tsize\frac{abcd}{(a,b,c,d)^4}\ \text{ is even}\\
{}&\text{if }\ \left(\tsize\frac{abcd}{(a,b,c,d)^4}, 6\right)\ne 1.
\endaligned
$$
Moreover
$$
\vth_a(\vth_{3/2})_b
\in J_{1,\frac{1}{2}(a^2+3b^2)}^{c}(v_\eta^4\times v_H^{a+b})
$$
if $(\frac a{(a,b)}, 3)=1$ or
$(\frac a{(a,b)},2)=2$ or $(\frac b{(a,b)},6)\ne 1$.
\endproclaim
\demo{Proof} The essential point is to prove that we have cusp forms.
For this we consider the norm of the indices of the non-zero Fourier
coefficients of the Jacobi forms written above.  For example, for the
third function we have
$$
\vartheta(\tau ,az)\vartheta(\tau ,bz)\vth(\tau,cz)=
\sum_{n,l,m\in \bz}
\biggl(\dsize\frac{-4}{n}\biggr)\biggl(\dsize\frac{-4}{m}\biggr)
\biggl(\dsize\frac{-4}{l}\biggr)
q^{\frac 1{8}({n^2+m^2+l^2})} r^{\frac 1{2}({an+bm+cl})}.
$$
Thus the norm of the index of $f(N,L)$ is given by
$$
2(a^2+b^2+c^2)N-L^2=
\frac 1{4}\bigl((bn-am)^2+(cn-al)^2+(cm-bl)^2\bigr).
$$
For $a$, $b$, $c$ satisfying the condition of the lemma the last sum
can only be zero if at least one of the three indices $n$, $m$ or $l$
is even.  For such $(n,m,l)$ the Fourier coefficient is $0$.  We have
the same formula for the norm of the indices of the fourth function.
In that case the norm is zero only if one of $n$, $m$ or $l$ has a
common divisor with 6.  Similarly, in the first two cases the norm is
a some of two squares, for the  fifth and sixth Jacobi form it is a
sum of four squares.
For the last function we have two different
Kronecker symbols in the product and the result follows
by similar arguments.

\hfill\hfill\qed
\enddemo

Using Jacobi forms of weight 1 from Lemma 1.2 and the lifting
construction gives us five series of $\Gamma_t$-cusp forms of weight $1$.

\proclaim{Lemma 1.3} The following cusp forms of weight $1$ exist:
$$
\spreadlines{2\jot}
\align
\al(\eta\vth_a)&\in \frak N_1(\Gamma_{3a^2}^+, \chi_{6,-}),
\qquad\quad
\al(\eta(\vth_{3/2})_a)\in \frak N_1(\Gamma_{18a^2}^+, \chi_{12,-})\\
\al(\vth_a\vth_b)&\in
\frak N_1(\Gamma_{2(a^2+b^2)}^+, \chi_{4, -})
\qquad
 \tsize\frac{ab}{(a,b)^2}\ \text{ is even}\\
\al((\vth_{3/2})_a(\vth_{3/2})_b)&
\in \frak N_1(\Gamma_{18(a^2+b^2)}^+, \chi_{12, -}),
\quad
\left(\tsize\frac{ab}{(a,b)^2}, 6\right)\ne 1\\
\al(\vth_a(\vth_{3/2})_b)&
\in \frak N_1(\Gamma_{3(a^2+3b^2)}^+, \chi_{6,-}),\quad
(\tsize \frac a{(a,b)}, 3)=1 \vee
(\tsize \frac a{(a,b)},2)=2 \vee
(\tsize \frac b{(a,b)},6)\ne 1.
\endalign
$$

\endproclaim

As we shall see later the cubes of the cusp forms of Lemma 1.3 define
canonical differential forms on the corresponding modular threefolds.
It is preferable to have cusp forms of weight $3$ with a character of
the smallest possible order. For example, the cubes of the first and
the fourth forms from Lemma \thetag{1.3} have a character of order $2$.  To
find differential forms on Siegel modular threefolds is one of the
main aims of this section.

Another way to construct new Jacobi forms of half-integral index is to
use the differential operators of Eichler-Zagier type.  For example if
$$
\phi_1\in J_{k_1,m_1}(v_\eta^{d_1}\times v_H^{\varepsilon_1}),
\qquad
\phi_2\in J_{k_2,m_2}(v_\eta^{d_2}\times v_H^{\varepsilon_2})
$$
are two Jacobi forms of integral or half-integral indices, where
$\varepsilon_i=0$ or $1$, then one can define the Jacobi form
$$
[\phi_1,\phi_2]=
\frac{1}{2\pi i}(m_2\pd{\phi_1}{z}\phi_2-m_1\phi_1\pd{\phi_2}{z})
\in J^{c}_{k_1+k_2+1,m_1+m_2}(v_\eta^{d_1+d_2}\times
v_H^{\varepsilon_1+\varepsilon_2}).
\tag{1.12}
$$
(See \cite{EZ, Theorem 9.5} and \cite{GN, Lemma 1.23}).  For arbitrary
integers $a$ and $b$ we obtain in this way the following Jacobi cusp
forms
$$
\gather
\phi_{2, \frac 1{2}(a^2+b^2)}(\tau,z)=
\frac{4}{ab}[\vth_a,\vth_b]=\\
\vspace{2\jot}
\sum_{m,n\in \bz}
(bm-an)\left(\frac{-4}m\right)\left(\frac{-4}n\right)\, q^{\frac
1{8}(m^2+n^2)}r^{\frac 1{2}(am+bn)}
\in J_{2, \frac1{2}(a^2+b^2)}^{c}(v_{\eta}^6\times v_H^{a+b}),
\tag{1.13}\\
\vspace{2\jot}
\phi_{3, \frac 1{2}(a^2+b^2+c^2+d^2)}(\tau,z)=
[\vth_a\vth_b,\,\vth_c\vth_d]
\in J_{3, \frac1{2}(a^2+b^2+c^2+d^2)}^{c}
(v_{\eta}^{12}\times v_H^{\Sigma})
\endgather
$$
where in the last formula $\Sigma:=a+b+c+d$. These Jacobi forms have
integral Fourier coefficients.

The Jacobi form $[\vth_a\vth_b,\,\vth_c\vth_d]$ is a cusp form of
weight $3$ with non-trivial binary character of the full Jacobi group.
We can define three different types of such Jacobi cusp forms which
give us altogether six series.  The first three series are
$$
\align
\eta^3\vth_a\vth_b\vth_c
\quad&\text{of index } \ \frac 1{2} ({a^2+b^2+c^2}),
\tag{I}\\
[\vartheta_a,\,\vartheta_b]\vth_c\vth_d
\quad&\text{of index } \  \frac 1{2}({a^2+b^2+c^2+d^2}),
\quad (a\ne b) \tag{II-a}\\
[\vth_a\vth_b,\,\vth_c\vth_d]
\quad&\text{of index } \  \frac 1{2}({a^2+b^2+c^2+d^2}).
\tag{II-b}
\endalign
$$
In the III-series we use the quintiple product:
$$
\align
\kern-30pt\eta^2(\vth_{3/2})_a\vth_{b}\vth_{c}\vth_{d}
\quad&\text{of index }\  \frac 1{2} ({3a^2+b^2+c^2+d^2}),
\tag{III-a}\\
\kern-30pt\eta(\vth_{3/2})_{a_1}(\vth_{3/2})_{a_2}
\vth_{b_1}\vth_{b_2}\vth_{b_3}\quad&\text{of index }\
\frac 1{2} ({3(a_1^2+a_2^2)+b_1^2+b_2^2+b_3^2}),
\tag{III-b}\\
\kern-30pt(\vth_{3/2})_{a_1}(\vth_{3/2})_{a_2}(\vth_{3/2})_{a_3}
\vth_{b_1}\vth_{b_2}\vth_{b_3}\quad&\text{of index }\
\frac 1{2} ({3(a_1^2+a_2^2+a_3^2)+b_1^2+b_2^2+b_3^2}).
\tag{III-c}
\endalign
$$

The III-series gives us new Jacobi forms if $a\ne b$, $c$, $d$ or
$a_i\ne b_j$.  The character of all Jacobi cusp forms of weight $3$
constructed in this way is $v_\eta^{12}\times v_H^{\Sigma}$, where
${\Sigma}$ is the sum of the corresponding indices $a$, $b$, $c$, $d$.

In \cite{G1} the first author constructed cusp forms of weight $3$
with respect to the paramodular group $\gt$ with trivial character for
all $t$ except
$$
t=1,\ 2,\,\dots,\,12,\ 14,\ 15,\ 16,\ 18,\ 20,\ 24,\ 30,\ 36.
\tag{1.14}
$$
We call these $t$ {\it exceptional}.  As a corollary it follows that
the geometric genus of any smooth model $\widetilde{\Cal A}_t$ of
${\Cal A}_t$ is positive if $t$ is not exceptional (see below). Here
we shall show the existence of many weight $3$ cusp forms with a
character for the exceptional values of $t$. Again we shall discuss
the geometric relevance of this below.

\proclaim{Theorem 1.4}
\roster
\item
Let $t$ be one of the exceptional polarizations and assume $t\not =1$,
$2$, $4$, $5$, $8$, $16$. Then there exists a cusp form of weight $3$
with respect to $\Gamma_t$ with the character $\chi_{2}$ of order $2$.
\smallskip
\item
For $t=8$, $16$ there exists a cusp form of weight $3$ with respect to
$\Gamma_t$ with the character $\chi_{4}$ of order $4$.
\smallskip
\item
Let $t\equiv 0\ \hbox{\rm mod}\ 3$ be exceptional and $t\ne 3,\,9$.
Then there exists a cusp form of weight $3$ with respect to $\Gamma_t$
with the character $\chi_{3}$ of order $3$.
\endroster
\endproclaim

\demo{Proof}
Let $\phi\in J_{3, \frac{1}2 t}^{c}(v_\eta^{12}\times v_H^{t})$, where
$t\in \bn$.  Then we can construct the lifting \thetag{1.9}
$$
\al(\phi)\in \frak N_{3}(\Gamma_t^+, \chi_{2, -}).
$$

Below we give a list of Jacobi forms of weight $3$ from the series
I--III for the exceptional polarizations \thetag{1.14}.  This will
prove the theorem.

We start the list of canonical cusp forms for the exceptional
polarizations with the cubes of cusp forms of weight $1$:
$$
\bigr(\al(\eta\vth)\bigl)^3=\Delta_1^3\in
\frak N_{3}(\Gamma_{3}^+, \chi_{2,-})
\tag{t=3}
$$
and
$$
\bigl(\al(\eta\vth_2)(Z)\bigr)^3
=F^{(12)}(Z)=\Delta_1( \left(\smallmatrix
\tau&2z\\2z&4\omega
\endsmallmatrix\right))^3\in
\frak N_{3}(\Gamma_{12}^+, \chi_{2, -}).
\tag{t=12}
$$
For $t=12$ there exists another cusp form, namely
$$
\al(\vth^2[\vth,\,\vth_3])\in
\frak N_{3}(\Gamma_{12}^+, \chi_{2, -})
\tag{t=12}
$$
(see the series II-a above).  To see that this is a different form one
can calculate the Fourier coefficients of the Jacobi forms we lift.

For $t=6$, $7$, $9$, $10$, $11$, $14$ and $24$ we obtain
$$
\spreadlines{2\jot}
\align
\al(\eta^3\vth^2\vth_2)&\in
\frak N_{3}(\Gamma_{6}^+, \chi_{2, -}),
\tag{t=6}\\
\al(\vth^2[\vth,\,\vth_2])&\in
\frak N_{3}(\Gamma_{7}^+, \chi_{2, -}),
\tag{t=7}\\
\al(\eta^3\vth\vth_2^2)&\in
\frak N_{3}(\Gamma_{9}^+, \chi_{2, -}),
\tag{t=9}\\
\al(\vth\vth_2[\vth,\,\vth_2])&\in
\frak N_{3}(\Gamma_{10}^+, \chi_{2, -}),
\tag{t=10}\\
\al(\eta^3\vth^2\vth_3)&\in
\frak N_{3}(\Gamma_{11}^+, \chi_{2, -}),
\tag{t=11}\\
\al(\eta^3\vth\vth_2\vth_3)&\in
\frak N_{3}(\Gamma_{14}^+, \chi_{2, -}),
\tag{t=14}\\
\al(\eta^3\vth^2\vth_2^2)&\in
\frak N_{3}(\Gamma_{24}^+, \chi_{2, -}).
\tag{t=24}
\endalign
$$
For the remaining values $t=12$, $15$, $18$, $20$, $30$, $36$ the
Jacobi forms of the series I--III provide us with several
$\Gamma_t$-cusp forms of weight $3$.  For $t=15$ we have five cusp
forms which are liftings of
$$
\vth\vth_3[\vth,\,\vth_2],\ \ \vth\vth_2[\vth,\,\vth_3],\ \
\vth^2[\vth_2,\,\vth_3], \ \ \eta^2\vth_{3/2}\vth_2^3,\  \
\eta^2(\vth_{3/2})_2\vth^3.
\tag{t=15}
$$
The liftings of these Jacobi forms generate a $4$-dimensional subspace
in $\frak N_{3}(\Gamma_{15}^+, \chi_{2, -})$.  We remark that the
first three Jacobi forms are obtained as ``permutations" of
theta-constants.

For $t=18$ the series I--III give us again five Jacobi cusp forms
$$
\eta^3\vth^2\vth_4,\  \
\eta(\vth_{3/2})^2\vth_2^3,\ \
\vth_2\vth_3[\vth,\,\vth_2]
\tag{t=18}
$$
(one can add two ``permutations" in the third Jacobi form).  Finally
 we get nine Jacobi cusp forms for $t=30$ with character
 $v_\eta^{12}\times 1_H$
$$
\eta^3\vth\vth_2\vth_5,\  \
\eta^2(\vth_{3/2})_2\vth^2\vth_4,\ \
\eta^2(\vth_{3/2})_3\vth^3,\ \
\eta(\vth_{3/2})^2\vth_2^2\vth_4,\ \
\vth_3\vth_4[\vth,\,\vth_2]
\tag{t=30}
$$
(one can add four ``permutations" of the fifth Jacobi form) and six
Jacobi forms for $t=36$ with the same character $v_\eta^{12}\times
1_H$
$$
\eta^3\vth_2\vth_4^2,\  \
\eta^2(\vth_{3/2})_3\vth\vth_2^2,\ \
\eta^2(\vth_{3/2})_3\vth^3,\ \
\vth\vth_5[\vth,\,\vth_3]
\tag{t=36}
$$
(plus two ``permutations" in the last Jacobi form).

A cusp form of weight $3$ for $\Gamma_8^+$ with the character
$\chi_{4,-}^{(8)}$ of order $4$ can be obtained, for example, by the
lifting of a Jacobi form of type II-a
$$
\al(\eta^2[\vth,\,\vth_{3/2}])
\in \frak N_3(\Gamma_8^+, \chi_{4, -}).
$$
To construct a similar cusp form for $\Gamma_{16}^+$ we consider the
Jacobi form
$$
\eta(\tau)^6\phi_{0,4}(\tau,z)=
\eta(\tau)^6\frac {\vth(\tau, 3z)}{\vth(\tau,z)}.
$$
The function $\phi_{0,4}(\tau,z)$ is a weak Jacobi form of weight $0$
and index $4$. One can prove that $\eta^2\phi_{0,4}(\tau,z)$ is a cusp
form (see \cite{GN, Lemma 1.21}). Thus
$$
\al(\eta^6\frac {\vth_3}{\vth})
\in \frak N_3(\Gamma_{16}^+, \chi_{4, -}).
$$

To prove the third claim of Theorem 1.4 we give a list of Jacobi cusp
forms of weight $3$ with the character $v_\eta^8\times1_H$:
$$
\aligned
\eta^5\vth_2\ \ &(t=6)\\
\eta^5\vth_4\ \ &(t=24)\\
{}
\endaligned
\qquad
\aligned
\eta[\vth,\,\vth_2]\vth_{3/2}\ \ &(t=12)\\
\eta^4\vth_3\vth_{3/2}\ \ &(t=18)\\
\eta^5{\vth_3\vth_4}/{\vth}\ \ &(t=36).
\endaligned
\qquad
\aligned
\eta^3\vth_2\vth_{3/2}^2\ \ &(t=15)\\
\eta^2[\vth_2,\,\vth_4]\ \ &(t=30)\\
{}
\endaligned
$$
The lifted forms are cusp forms of weight $3$ with the character
$\chi_{3,-}$.

\hfill\hfill\qed
\enddemo

At this point it is natural to discuss the geometric implications of
our results. Weight $3$ cusp forms are closely related to canonical
differential forms on smooth models of the corresponding modular
variety. If $F$ is a cusp form of weight $3$ with respect to a group
$\Gamma$, then $\omega_F=F(Z)dZ$ is a holomorphic $3$-form on the
space $\Cal A_\Gamma^o=(\Gamma\setminus \bh_2)^o$, where ${}^o$ means
that we consider the threefold outside the branch locus of the natural
projection from $\bh_2$ to $\Cal A_\Gamma$. A very useful extension
theorem due to E. Freitag implies that such a form can be extended to
any smooth model of $\Cal A_\Gamma$. To be more precise, let $\Gamma$
be an arbitrary subgroup of $Sp_4(\br)$, which contains a principal
congruence subgroup $\Gamma_1(q)\subset Sp_4(\bz)$ of some level
$q$. We then have the following
\proclaim{Criterion}(Freitag) An element
$\omega_F=F(Z)dZ \in H^0(\Cal A_\Gamma^o,\, \Omega_3(\Cal
A_\Gamma^o))$ can be extended to a canonical differential form on a non-singular
model $\
\widetilde{\Cal A}_\Gamma\ $ of a compactification of $\Cal A_\Gamma$
if and only if the differential form $\omega_F$ is square integrable.
\endproclaim
\demo{Proof}See \cite{F}, Hilfsatz 3.2.1.
\enddemo

It is well known that a $\Gamma$-invariant differential form
$\omega_F=F(Z) dZ$ is square-integrable if and only if $F$ is a cusp
form of weight $3$ with respect to the group $\Gamma$.  Thus we have
the following identity for the geometric genus of the variety
$\widetilde{\Cal A}_\Gamma$:
$$
p_g(\widetilde{\Cal A}_\Gamma)=h^{3,0}(\widetilde{\Cal A}_\Gamma)
=\hbox{dim}_\Bbb C \,\frak N_3(\Gamma).
\tag{1.15}
$$
We also remark at this point that $\Gamma_t$-cusp forms of weight $2$
can be very useful when one wants to prove that some modular
threefolds are of general type (see \cite{GH1}, \cite{GS} and
\cite{S}).  All
these facts explain our interest in cusp forms of small weight $k$
($k\le 3$). To reformulate our above results in geometric terms we
introduce the following
\definition{Definition 1.5}Let $\Gamma$ be a subgroup of $Sp_4(\br)$
which contains a principal congruence-subgroup.  For any character
$\chi_\Gamma:\Gamma\to \bc^*$ we define the threefold
$$
\Cal A(\chi_\Gamma)=\hbox{ker} (\chi_\Gamma)\setminus \bh_2.
$$
The covering $
\ \Cal A(\chi_\Gamma)\to \Cal A_\Gamma=\Gamma\setminus \bh_2\
$ is galois with a finite abelian Galois group.
\enddefinition
\medskip

\proclaim{Corollary 1.6}Let $t$ be one of the exceptional
polarizations (see \thetag{1.14}).
\roster
\item
If $t\ne 1,\,2,\,4,\,5,\,8,\,16$, then the modular double covering
$$
\Cal A(\chi_{2})\overset{2:1}\to\longrightarrow \Cal A_t
$$
of the moduli space $\Cal A_t$ of abelian surfaces with polarization
of type $(1,t)$ has positive geometric genus, and in particular the Kodaira
dimension of ${\Cal A}(\chi_2)$ is not negative. Moreover, for
$t=12$,
$15$, $18$, $20$, $30$ and $36$ the Kodaira dimension of $\Cal A_t$ is
positive.
\smallskip
\item
If $t=6,\,12,\,15,\,18,\,24,\,30,\,36$, then the threefold $
\Cal A_t(\chi_{3})\overset{3:1}\to\longrightarrow \Cal A_t
$ has positive geometric genus.
\smallskip
\item
If $t=8$ or $16$, then the covering $
\Cal A_t(\chi_{4})\overset{4:1}\to\longrightarrow \Cal A_t
$ has positive geometric genus.
\endroster
In particular all these modular varieties are not unirational.
\endproclaim
\remark{Remark 1.7}If one lifts a Jacobi form which contains a factor
$\vth$ or $\vth_{3/2}$, then one has some information about
the divisor of the lifted form (see \cite{GN, Lemma 1.16}).  Assume,
e.g. that $\vth(\tau, az)$ is a factor of a Jacobi form $\phi$, then
$$
\dv_{\Cal A_t^+}(\al(\phi))\supset
\sum_{d|a} H_{d^2}(d).
$$
Here we denote by
$$
H_D(b)=\pi_t^+(\{Z=\left(\smallmatrix \tau&z\\z&\omega
\endsmallmatrix\right)\in \bh_2\,|\,
a\tau+bz+t\omega=0\})
$$
the Humbert modular surface of discriminant $D=b^2-4at$ in $\Cal
A_t^+$ where
$$
\pi_t^+: \bh_2\to \Cal A_t^+=\Gamma_t^+\setminus \bh_2
$$
is the natural projection (see \cite{GH2}).  For example, we obtain
$$
\dv_{\Cal A_{15}^+}(\al(\eta^2\vth_{3/2}\vth_2^3))
\supset 3H_1+4H_4(2),\ \
\dv_{\Cal A_{15}^+}(\al(\vth\vth_3[\vth,\,\vth_2]))
\supset 3H_1+H_9(3).
$$
The other cusp forms of weight $3$ constructed above can be treated
similarly.

The modular form $\Delta_1=\al(\eta\vth)$ was studied in \cite{GN},
where its was proved that its divisor in $\Cal A_3^+$ is exactly equal
to the Humbert modular surface with discriminant $1$.  Thus
$$
\dv_{\Cal A_{3}^+}(\al (\eta\vth)^3)=3H_1, \qquad
\dv_{\Cal A_{12}^+}(\al (\eta\vth_2)^3)=3H_1+3H_4(2).
$$
On the other hand Brasch \cite{B} has determined the branch locus of
the map
$$
\pi_t: \bh_2\to \Cal A_t=\Gamma_t\setminus \bh_2
$$
as well as the singularities of a toroidal compactification of $\Cal
A_t$ (at least in the prime number case).  From his work it is
possible to obtain the same information about coverings of $\Cal A_t$.
Altogether this makes it possible to write down an explicit effective
canonical divisor on a suitable smooth model.  This information is
very useful if one wants to study the geometry of these modular
varieties from the point of view of Mori theory. We also hope to find
 examples of modular varieties which are Calabi--Yau
varieties. We hope to return to this in the future.
\endremark

\head
\S 2. Commutator subgroups
\endhead

The aim of this section is to determine the commutator subgroup of
$\Gamma_t$ and $\Gamma_t^+$. In the case of $\Gamma_1$ this is a
classical result due to Reiner \cite{R} and Maa\ss\ \cite{Ma}.  We
shall also discuss some consequences of our result.
\proclaim{Theorem 2.1}
For any integer $t\ge 1$ let $t_1=(t,12)$ and $t_2=(2t,12)$. If
$\Gamma_t'$ is the commutator subgroup of $\Gamma_t$,
resp. $(\Gamma_t^+)'$ is the commutator subgroup of $ \Gamma_t^+$,
then the following holds:
\roster
\item"{\rm (1)}"
$\Gamma_t/\Gamma_t'\cong \bz/t_1\times \bz/t_2$
\smallskip
\item"{(\rm 2)}"
$\Gamma_t^+/(\Gamma_t^+)'\cong \bz/2 \times \bz/t_2.$
\endroster
\endproclaim
\demo{Proof}
We shall proceed in several steps.
\smallskip

{\bf Step 1:} Recall that we have the following two embeddings of
$SL_2(\bz)$ into the paramodular group $\Gamma_t$:
$$
i_{\infty}\pmatrix a & b\\ c & d
\endpmatrix=
\pmatrix
a & 0 & b & 0\\ 0 & 1 & 0 & 0\\ c & 0 & d & 0\\ 0 & 0 & 0 & 1
\endpmatrix, \qquad
j_{\infty}
\pmatrix
a' & b'\\ c' & d'
\endpmatrix=
\pmatrix
1 & 0 & 0 & 0\hphantom{{}^{-1}}\\ 0 & a' & 0 & b't^{-1}\\ 0 & 0 & 1 &
0\hphantom{{}^{-1}}\\ 0 & c't & 0 & d'\hphantom{{}^{-1}}
\endpmatrix.
$$
We consider the matrices
$$
A=i_{\infty}\pmatrix 1 & 1\\ 0 & 1
\endpmatrix=
\pmatrix
1 & 0 & 1 & 0\\ 0 & 1 & 0 & 0\\ 0 & 0 & 1 & 0\\ 0 & 0 & 0 & 1
\endpmatrix,
\
B=j_{\infty}\pmatrix 1 & 1\\ 0 & 1
\endpmatrix=
\pmatrix
1 & 0 & 0 & 0\hphantom{{}^{-1}}\\ 0 & 1 & 0 & t^{-1}\\ 0 & 0 & 1 &
0\hphantom{{}^{-1}}\\ 0 & 0 & 0 & 1\hphantom{{}^{-1}}
\endpmatrix.
\tag{2.1}
$$
It is well known \cite{Ma, p. 130} that the commutator $SL_2(\bz)'$ of
$SL_2(\bz)$ has index 12 in $SL_2(\bz)$ and that
$$
SL_2(\bz) =
\bigcup^{11}_{l=0}SL_2(\bz)'T^l
$$
\noindent
where $T=\left(\smallmatrix 1 & 1\\ 0 & 1
\endsmallmatrix\right)$.
Applying this to the two subgroups in $\Gamma_t$ isomorphic to
$SL_2(\bz)$ via $i_{\infty}$, resp. $j_{\infty}$ it follows that in
particular $A^{12}, B^{12} \in \Gamma_t'$.
\smallskip
As in \cite {Ma} we can construct elements in $\Gamma_t'$ in the
following way
$$
\pmatrix
V & 0\\ 0 & ^tV^{-1}
\endpmatrix
\pmatrix
{\bold 1}_2 & S\\ 0 & {\bold 1}_2
\endpmatrix
\pmatrix
V^{-1} & 0\\ 0 & ^tV
\endpmatrix
\pmatrix
{\bold 1}_2 & -S\\ 0 & {\bold 1}_2
\endpmatrix=
\pmatrix
{\bold 1}_2 & VS{^tV}-S\\ 0 & {\bold 1}_2
\endpmatrix
$$
where $V=\pmatrix a & tb\\ c & d
\endpmatrix
\in GL_2(\bz)
$ and $S=\pmatrix m & n\\ n & t^{-1}k
\endpmatrix$.
In particular, we obtain for $W=VS{^tV}-S$ the following matrices.
\roster
\item
$ W=\pmatrix 2tn+tk & k\\ k & 0
\endpmatrix\qquad \hbox { for }
V=
\pmatrix
1 & t\\ 0& 1
\endpmatrix
$
\smallskip
\item
$ W=\pmatrix 0 & -2n\\ -2n & 0
\endpmatrix\qquad \hbox { for }
V=\pmatrix 1 & 0\\ 0& -1
\endpmatrix
$
\smallskip
\item
$ W=\pmatrix 0 & m\\ m & m+2n
\endpmatrix\qquad
\hbox { for }
V=\pmatrix 1 & 0\\ 1& -1
\endpmatrix.
$
\endroster
Together with $A^{12}, B^{12}\in \Gamma_t'$ we find that $A^{t_2},
B^{t_2} \in \Gamma_t'$.

\noindent {\bf Step 2:}
We consider the following subgroup of $\Gamma_t$:
$$
G=\bigcup^{t_2-1}_{l,m=0} \Gamma_t' A^lB^m.
$$

\proclaim{Claim}
$G=\Gamma_t$.
\endproclaim

To prove this claim we consider the matrix
$$
C=\pmatrix 1 & 0 & 0 & \hphantom{-}0\\ 1 & 1 & 0 & \hphantom{-}0\\ 0 &
0 & 1 & -1\\ 0 & 0 & 0 & \hphantom{-}1
\endpmatrix.
$$
Our first aim is to show that $C$ is contained in $G$.  Indeed,
consider
$$
L=\pmatrix 1 & 0 & -1 & 0\\ 0 & 1 & \hphantom{-}0 & 1\\ 0 & 0 &
\hphantom{-}1 & 0\\ 0 & 0 & \hphantom{-}0 & 1
\endpmatrix
=A^{-1} B^t \in G, \quad
M=\pmatrix
\hphantom{-}1 & 0 & 0 & 0\\
\hphantom{-}0 & 1 & 0 & 0\\
-1 & 0 & 1 & 0\\
\hphantom{-}0 & 0 & 0 & 1
\endpmatrix
\in \Gamma'_t A^{-1} \subset G
$$
(the latter follows from the corresponding statement in
$SL_2(\bz)$ via the inclusion $i_{\infty}$).
$$
X=\pmatrix
\hphantom{-}0 & 0 & 1 & 1\\
 -1 & 1 & 1 & 1\\ -1 & 0 & 2 & 1\\
\hphantom{-}0 & 0 & 0 & 1
\endpmatrix \in \Gamma_t.
$$
By a straightforward calculation
$$
C=L A X A^{-1} X^{-1} M\in G.
$$
\noindent
Recall from \cite{G1, Lemma 2.2} that
$\Gamma_{t, \infty}=\Gamma_t\cap \Gamma_\infty(\bq)$ and the
element
$$
J_t =\pmatrix
\hphantom{-}0 & \hphantom{-}0 & 1 & 0\hphantom{{}^{-1}}\\
\hphantom{-}0 & \hphantom{-}0 & 0 & t^{-1}\\
-1 &\hphantom{-} 0 & 0 & 0\hphantom{{}^{-1}}\\
\hphantom{-}0 & -t            & 0 & 0\hphantom{{}^{-1}}
\endpmatrix
$$
generate $\Gamma_t$.  Hence it is enough to show that
$\Gamma_{t, \infty}$ and $J_t$ are contained in $G$.  The assertion
about $J_t$ follows since both copies of $SL_2(\bz)$ which are
contained in $\Gamma_t$ are already contained in $G$.  Now consider
$g\in\Gamma_{t, \infty}$.  Again using the two copies of $SL_2(\bz)$ we
can assume that
$$
g=\pmatrix 1 & 0 & 0 & \hphantom{{-}}k\\ r & 1 & k& \hphantom{{-}}0\\
0 & 0 & 1 & -r\\ 0 & 0 & 0 & \hphantom{{-}}1
\endpmatrix.
$$
But then
$$
g C^{-r} B^{-rkt} A^{kt}=
\pmatrix
1 & 0 & k & k\\ 0 & 1 & k & 0\\ 0 & 0 & 1 & 0 \\ 0 & 0 & 0 & 1
\endpmatrix \in \Gamma'_t
$$
by step 1 and this gives the claim.
\smallskip
\noindent {\bf Step 3:} We shall now prove the theorem for $\Gamma_t^+$.
Recall from Lemma 1.1 that there exist characters
$$
\chi_{t_{2,\pm}}:\Gamma_t^+\rightarrow \{\root{t_2}\of{1}\}
$$
of order $t_2$ with $\chi_{t_{2,\pm}}|_{SL_2(\bz)}=v^{24/t_2}_{\eta}$,
$\chi_{t_{2,\pm}}|_{H({\bz})}=v_{H}^{\varepsilon}$ with $\varepsilon=0$ or
$1$ depending on $t$ and $\chi_{t_2,\pm}(V_t)=\pm 1$. Since
$\Gamma_t$ is a normal subgroup of $\Gamma_t^+$ of index $2$, there also exists
a character of order two:
$$
\chi'_2:\Gamma_t^+\rightarrow\{\pm 1\}
$$
with ${\chi'_2}|_{\Gamma_t}=1$. Clearly $\chi_{t_{2,-}}=\chi_{t_{2,+}}
\chi'_2$. Together
$\chi_{t_{2,+}}$ and
$\chi'_2$ define a surjective map
$$
{\bar \chi}=(\chi'_2, \chi_{t_{2,+}}):
\Gamma_t^+\rightarrow \bz/2\times \bz/t_2.
$$
On the other hand $B=V_t AV_t^{-1}$, i.e. $\chi(A)=\chi(B)$ for every
character $\chi$ of $\Gamma_t^+$.  Together with step 2 this gives the
claim.
\smallskip
\noindent {\bf Step 4:} It remains to prove the theorem for $\Gamma_t$.
In view of step 2 this will be a consequence of the following facts:
\roster
\item
There exists a character $\chi_{t_2}$ of $\Gamma_t$ of order $t_2$ with
$\chi_{t_2}(A)=\chi_{t_2}(B)=e^{2\pi i/t_2}$.
\smallskip
\item
There exist characters $\chi_{t_1}$ and $ \chi_{t_1}'$ of $\Gamma_t$ with
$\chi_{t_1}(A)=e^{2\pi i/t_1}, \chi_{t_1}(B)=1$ and $\chi'_{t_1}(A)=1,
\chi'_{t_1}(B)=e^{2\pi i/t_1}$.
\smallskip
\item
For every character $\chi$ of $\Gamma_t$ the equality
$\chi(A^{t_1})=\chi(B^{t_1})$ holds.
\endroster
Assertion (1) is clear, since we can take the character $\chi_{t_{2,\pm}}$
of order $t_2$ on $\Gamma_t^+$ and restrict to $\Gamma_t$.  To prove
(2) it is slightly more natural to work with the group ${\widetilde
\Gamma}_t=Sp(W_t,\bz)$.  Recall that if $g\in{\widetilde\Gamma}_t$,
then
$$
g=
\pmatrix
a & t{\ast} & b & t{\ast}\\
\ast  & \ast & \ast & \ast\\
c & t{\ast} & d & t\ast \\
\ast  & \ast & \ast & \ast
\endpmatrix
\tag{2.2}
$$
where $a,\ldots , d$ and $\ast$ are integer entries (because of a bad
line break).  Since the matrix $g$ is symplectic it follows that $ad-bc \equiv
1\mod t$, i.e.
$$
\pmatrix
a & b\\ c & d
\endpmatrix
\ \hbox{mod}\  t \in SL_2( \bz/t).
$$
On the other hand $v_{\eta}$ defines a character
$$
\aligned
v_{\eta}^{24/t_1}: &\ SL_2( \bz/t)\\
\vspace{2\jot}
{}&\ \pmatrix a & b\\ c & d
\endpmatrix \mod t
\endaligned
\aligned
{}& \longrightarrow
\{ \root{t_1}\of{1}\}\\
\vspace{2\jot}
{} & \longmapsto v_{\eta}^{24/t_1} (\pmatrix a & b\\ c & d
\endpmatrix).
\endaligned
$$
For every element $g\in {\widetilde\Gamma}_t$ we set
$$
\chi_{t_1}(g):=v_{\eta}^{{24}/{t_1}}\left(
\pmatrix
a & b\\ c & d
\endpmatrix \mod t \right).
$$
It follows easily from \thetag{2.2} that this is indeed a character.
The character $\chi'_{t_1}$ can then be defined by
$\chi'_{t_1}(g)=\chi_{t_1}(V_tg V_t^{-1})$. To prove (3) we recall
from step 1 (take $n=0, k=-m=1)$ that
$$
A^{t} B^{-t} =
\pmatrix
1 & 0 & t & \hphantom{{-}}0\\ 0 & 1 & 0 & -1\\ 0 & 0 & 1 &
\hphantom{{-}}0 \\ 0 & 0 & 0 & \hphantom{{-}}1
\endpmatrix \in \Gamma'_t.
$$
But then also $A^{t_1} B^{-t_1}\in \Gamma'_t$ and hence
$\chi(A^{t_1})=\chi(B^{t_1})=1$ for every character $\chi$ of
$\Gamma_t$.

\hfill\hfill\qed
\enddemo

\noindent{\bf Notation 2.2.} The proof of Theorem 2.1 shows that a character
$\chi:\Gamma_t\rightarrow \bc^{\ast}$ is determined by its values on $A$ and
$B$. These values can be any root of order $t_2$ provided
$\chi(A^{t_1})=\chi(B^{t_1})$. We denote by
$\chi_{a,b}:\Gamma_t\rightarrow\bc^{\ast}$  ($1\le a, b\le t_2$ and
$(a-b) \equiv t_2/t_1$)
the unique character of $\Gamma_t$ given by
$\chi_{a,b}(A)=e^{2\pi i a/t_2}$ and $
\chi_{a,b}(B)=e^{2\pi i b/t_2}$. Comparing this with the definition of $\chi_Q$
from section 1 shows that $\chi_Q=\chi_{a,a}$ if and only if $aQ=t_2$.

Similarly a character $\chi:\Gamma^+_t\rightarrow \bc^{\ast}$ is determined by
its values on $A$ and $V_t$. We denote by $\chi_{a,a,\pm}$ the unique character
of $\Gamma^+_t$ such that $\chi_{a,a,\pm}(A)=e^{2\pi i a/t_2}$ and
$\chi_{a,a,\pm}(V_t)=\pm1$. Then $\chi_{a,a,\pm}|\Gamma_t=\chi_{a,a}$.
\medskip

The classical construction of the character $\chi_2$ of order $2$ for
$Sp_4 (\bz)$ can be easily generalized to $\Gamma_t$ for $t$ odd.
Again it is more convenient to work with the group
${\widetilde\Gamma}_t=Sp (W_t,\bz)$.  Let
$$
Sp^{(2)}(W_t,\bz)=
\{g\in Sp (W_t,\bz),\  g \equiv{\bold 1}\ \hbox{mod}\  2\}.
$$
Then one can show as in \cite{Ig1, Lemma V.2.5} that there is an exact
sequence
$$
1\rightarrow Sp^{(2)}(W_t,\bz)\rightarrow Sp(W_t,\bz)\rightarrow
Sp_4(\bz/2)\rightarrow 1.
$$
The group $Sp_4(\bz/2)$ is isomorphic to the symmetric group $S_6$ and
hence has a character of order $2$. This defines a character of
order $2$ for $Sp(W_t,\bz)$ and hence also for $\Gamma_t$.  The moduli
space defined by $Sp^{(2)}(W_t,\bz)$ is the moduli space of
$(1,t)$-polarized abelian surfaces with a level-$2$ structure.
$$
\align
{\Cal A}^{(2)}_t {}&= {Sp}^{(2)}(W_t,\bz)\backslash \bh_2\\
{}&=\{(A,H,\alpha);\ H \hbox { is a } (1,t)\hbox{-polarization, } \
\alpha \hbox{ is a level-$2$ structure}\}.
\endalign
$$
The group $Sp_4(\bz/2)\cong S_6$ acts transitively on the set of
level-$2$ structures of a fixed abelian surface.  The alternating
group $A_6$ has two orbits. We call these orbits {\it classes} of
level-$2$ structures.  So if $\ker(\chi_2)\subset \Gamma_t$ is the
kernel of the character of order $2$ constructed above, then
$$
\align
{\Cal A}(\chi_2)&=\ker(\chi_2) \backslash \bh_2\\ {}&=\{
(A,H,{[\alpha]});
\  H \hbox{ is a } (1,t)\hbox{-polarization, }
{\ [\alpha]} \hbox{ is a class of level-$2$ structures} \}.
\endalign
$$

We want to conclude this section with an application of Theorem 2.1 to
the Picard group of the moduli stack of (1,t)-polarized abelian
surfaces.  In his beautiful paper
\cite{Mu1} Mumford computes the Picard group of the moduli stack
${\Cal M}_{1,1}$ of elliptic curves.  His result is that $\hbox{
Pic}({\underline{\Cal M}}_{1,1})=\bz/12$ which is the quotient of
$SL_2(\bz)$ by its commutator subgroup.  Similarly the computation of
the commutator group of $\Gamma_t$ has an interpretation in terms of
the Picard group of the moduli stack ${\underline{\Cal A}}_t$ of
$(1,t)$-polarized abelian surfaces.

\proclaim{Proposition 2.3}
The Picard group $\operatorname{Pic}({\underline{\Cal A}_t})$ is
finitely generated of rank $\ge 1$. Moreover
\roster
\item"{\rm (i)}"
$\operatorname {Pic}({\underline{\Cal A}}_1)\cong \bz \times
\bz/2$
\smallskip
\item"{\rm (ii)}"
$\operatorname {Tor} \operatorname {Pic}({\underline{\Cal A}}_t)\cong
\bz/t_1\times
\bz/t_2$ \quad for  $t>1$.
\endroster
\endproclaim

\demo{Proof}
Let ${\Cal L}$ be the $\bq$-line bundle associated to modular forms of
weight $1$. We first remark that ${\Cal L}$ is a non-zero element of
$\operatorname {Pic}({\underline{\Cal A}}_t)\otimes \bq$.  In fact
${\Cal L}$ is non-zero whenever the genus $g\ge 2$.  Assume that
${\Cal L}^{\otimes k}= {\Cal O} $ for some $k$.  Then ${\Cal
L}^{\otimes kn}={\Cal O}$ for all $n\ge 1$.  But a trivializing
section for ${\Cal L}^{\otimes kn}$ would imply the existence of a
non-constant modular form without zeroes, contradicting Koecher's
principle for $g\ge 2$. In the principally polarized case it is known
(see e.g. \cite{Mu2, Part III}) that ${\Cal L}$ generates
$\operatorname {Pic}({\Cal A}_1)\otimes \bq$ and hence the same is
true for $\operatorname {Pic}({\underline{\Cal A}}_1)\otimes \bq$.
(This follows exactly as in \cite {C, Lemme 1.1}.)
For the rest of the proof we can now argue very much as in \cite{Mu1,
\S 7} and hence we omit full details.  First of all we have an exact
sequence
$$
\aligned
0 &\rightarrow H^0({\underline{\Cal A}}_t,\bz/n)
\rightarrow
H^0({\underline{\Cal A}}_t,{\Cal O}^{\ast})\\ {}&\rightarrow
H^1({\underline{\Cal A}}_t,\bz/n)
\rightarrow
\operatorname {Pic}({\underline{\Cal  A}}_t)
\endaligned
\aligned
{}&
\overset {n}\to{\rightarrow}
H^0({\underline{\Cal A}}_t, {\Cal O}^{\ast})\\ {}&\overset
{n}\to{\rightarrow}
\operatorname {Pic}({\underline{\Cal  A}}_t)
\rightarrow
H^1({\underline{\Cal A}}_t,\bz/n).
\endaligned
$$
Moreover
$$
H^1({\underline{\Cal A}}_t, \bz/n) = H^1(\Gamma_t, \bz/n)
=\operatorname{Hom}(\Gamma_t/\Gamma_t', \bz/n).
$$
By our computation of $\Gamma_t'$ it follows that
$H^1({\underline{\Cal A}}_t, \bz/p)=0$ for every prime $p$ with $(p,
2t)=1$.  As in \cite {Mu1, p.77} this implies that
$H^0({\underline{\Cal A}}_t,{\Cal O}^{\ast})=\bc^{\ast}$ and that
$\operatorname {Pic}({\underline{\Cal A}}_t)$ is finitely generated.
But now we have an exact sequence
$$
0\rightarrow H^1({\underline{\Cal A}}_t, \bz/n)\rightarrow
\operatorname
{Pic}({\underline{\Cal A}}_t) \overset {n}\to{\rightarrow}
\operatorname {Pic}({\underline{\Cal  A}}_t)
$$
\noindent
and since the $n$-torsion of $\operatorname{Pic }({\underline{\Cal
A}}_t)$ is just the kernel of the last map in this sequence, the
result follows from Theorem 2.1.

\hfill\hfill{\qed}
\enddemo

\remark{Remark}
Geometrically the 2-torsion element in $\operatorname
{Pic}({\underline{\Cal A}}_1)$ can be realized as follows: Let
$H_1=\pi_1\{\left(\smallmatrix \tau&0\\0&\omega\endsmallmatrix\right)\}
\subset {\Cal A}_1$
be the Humbert surface of discriminant $1$.
Then $H_1=\pi_1\{\Delta_5=0\}$ where $\Delta_5(Z)$ is the
product of the ten even theta-characteristics. The function
$\Delta_5(Z)$ is a modular form of weight 5 with a non-trivial
character of order $2$.  Hence $2H_1=10{\Cal L}$ in $\operatorname
{Pic}({\underline{\Cal A}}_1)$, but $H_1$ and $5{\Cal L}$ differ by a
non-zero 2-torsion element.
\endremark
\medskip

\head
\S 3. The geometric genus of the commutator covering of $\Cal A_t$
\endhead
\smallskip

In this section we consider modular coverings of the moduli space
$\Cal A_t$ with commutative covering groups.  The commutator subgroup
$\Gamma_t'\subset \Gamma_t$ defines the maximal abelian covering
$$
\Cal A_t^{com}=\Gamma_t'\setminus \bh_2 \to
\Gamma_t\setminus \bh_2=\Cal A_t.
\tag{3.1}
$$
We denote by $\widetilde{\Cal A}$ an arbitrary smooth compact model of
a threefold $\Cal A$.

\proclaim{Theorem 3.1}
\roster
\item
The maximal abelian covering $\widetilde{\Cal A}_t^{com}$ has
geometric genus $0$ if and only if $t=1,\ 2,\ 4,\ 5$.
\smallskip
\item The threefolds
${}\quad{\Cal A}_3^{com}\overset{18:1}\to \longrightarrow {\Cal A}_3\
$, $\ {\Cal A}_7^{com}\overset{2:1}\to \longrightarrow {\Cal A}_7\ $
and $\ {\Cal A}_6(\chi_{2})
\overset{2:1}\to \longrightarrow
{\Cal A}_6\quad$ have
\newline
geometric genus $1$: $\ h^{3,0}(\widetilde{\Cal A}_3^{com})=
h^{3,0}(\widetilde{\Cal A}_7^{com})=
h^{3,0}(\widetilde{\Cal A}_6(\chi_{2}))=1$.
\smallskip
\item
Let $\chi_3:\Gamma_6\to \bc^*$ be the character \thetag{1.10a--b}
of order $3$.
Then $h^{3,0}(\widetilde{\Cal A}_6(\chi_{3}))$=2.
\endroster
\endproclaim

We shall prove the theorem using the lifting construction \thetag{1.9}
together with a statement about the $q$-order of Jacobi forms which we
formulate below.  Let
$$
\phi_{k,m}(\tau,z)=
\sum\Sb n,l\in \bz\\
\vspace{0.5\jot}
4mn-l^2>0
\endSb
a(n,l)q^nr^l\in J_{k,m}^c
$$
be a Jacobi cusp form with trivial character.  The {\it $q$-order} of
$\phi_{k,m}(\tau,z)$ is the minimal $q$-power contained in its Fourier
expansion
$$
\hbox{ord}_q(\phi_{k,m})=\hbox{min}\,\{n\in \bn\,|\,a(n,l)\ne 0\}.
$$
\proclaim{Proposition 3.2}
Let $\phi_{2k,m}\in J_{2k,m}^c$ be a non-zero Jacobi cusp form of even
weight. Then
$$
\hbox{\rm ord}_q(\phi_{2k,m})\le \hbox{\rm min}\,
(\frac{3k-3+m}9, \ \frac{k+m}6).
\tag{3.2}
$$
Moreover equality in \thetag{3.2} holds if and only if
$\phi_{2k,m}=c\Delta^N\eta^{-6m}\vth^{2m}$ where $N> \frac m{4}$ is
the $q$-order, $2k=12N-2m$ is the weight, and $c\ne 0$ is a constant.
\endproclaim

\medskip

The estimation of this proposition works very well for applications
concerning modular forms with respect to $\Gamma_t$ for $t\le 8$.  For
example we easily obtain the following

\proclaim{Corollary 3.3}For $t\le 8$ the geometric genus
of the moduli space of abelian surfaces with a polarization of type
$(1,t)$ is zero:
$$
h^{3,0}(\widetilde{\Cal A}_t)=0.
$$
\endproclaim
\demo{Proof of the Corollary}
According to \thetag{1.15} we have to prove that there are no
$\Gamma_t$-cusp forms of weight $3$ for $t\le 8$.  Consider the
decomposition
$$
\frak N_k(\Gamma_t)=
\frak N_k^{+}(\Gamma_t)\oplus \frak N_k^-(\Gamma_t),
\tag{3.3}
$$
where the $\pm$-subspaces consist of $\Gamma_t$-modular forms, which
are invariant or
anti-invariant with respect to the involution $V_t$ (see \thetag{1.1})
$$
F(Z)=\pm F(V_t<Z>),
\qquad V_t<Z>=\pmatrix \omega/t &z\\z&t\tau
\endpmatrix.
$$
Because of the decomposition \thetag{3.3} we can restrict ourselves to
$\Gamma_t^+$-cusp forms.

Let $F_3^{(t)}\in \frak N_3^{\pm}(\Gamma_t)$.  We consider its
Fourier-Jacobi expansion
$$
F_3^{(t)}(Z)=\sum\Sb m\ge M>0\endSb f_{3, {mt}}(\tau,z)e^{2\pi i m
t\omega}.
$$
From $V_t$-invariance of $\bigl(F_3^{(t)}\bigr)^2$ and Proposition 3.2
follows that $ 2M\le \hbox{ord}_q(f_{3,Mt}^2)<\frac{6+2Mt}9 $.  For
$t\le 6$ we have $\frac{6+2Mt}9\le 2M$, thus there are no cusp forms
of weight $3$ for $t\le 6$.  For $t=7$ and $t=8$ the same inequality
holds for $M\ge 2$ and $M\ge 3$ respectively. We note that
$J_{3,7}=\{0\}$ and $J_{3,8}^c=J_{3,16}^c=\{0\}$ (see below).  This
finishes the proof.

\hfill\hfill\qed
\enddemo

Recall that dimension formulae for $J_{k,m}$ were found in \cite{EZ}
and \cite{SkZ}. We use these in the following variant
$$
\hbox{dim}_\Bbb C\, J_{k,m}^{c}
=\cases
\sum_{j=0}^{m}\{k+2j\}_{12}
-\lfloor\frac{j^2}{4m}\rfloor \qquad &k\text{ is even}\\
\sum_{j=1}^{m-1}\{k+2j-1\}_{12}
-\lfloor\frac{j^2}{4m}\rfloor \qquad &k\text{ is odd}
\endcases
\tag{3.4}
$$
with
$$
\{m\}_{12}=\cases\lfloor\frac{m}{12} \rfloor&\text{if }\ m\not\equiv 2
\mod 12\\
 \lfloor \frac{m}{12} \rfloor-1&\text{if }\ m\equiv 2 \mod 12.
                       \endcases
$$

\remark{Remark}The fact that
$\hbox{dim}\,\frak N_3^{+}(\Gamma_t)=0$ for
$t=2,\, 3, \, 5,\, 7$ was
proved in \cite{G2, \S 3.5} by a different method.  Gross and Popescu
(cf. \cite{GP}) have announced a result which, among other things,
implies that the moduli spaces $\Cal A_t$ are unirational for $t \le
12$. Their approach is via equations for abelian surfaces.
\endremark
\demo{Proof of Proposition 3.2}
Let $N=\hbox{ord}_q(\phi_{2k,m})$.  If $6N<k$, then the proposition is
proved.  Assume that $6N\ge k$ and consider the function
$\Delta(\tau)^{-N}\phi_{2k,m}(\tau, z)$, where
$\Delta(\tau)=\eta(\tau)^{24}$.  According to the definition of $N$
the function $\Delta^{-N}\phi_{2k,m}$ is a weak Jacobi form of weight
$2(k-6N)$ (see \thetag {1.2}).  By \cite{EZ, Theorem 9.3} the ring
$J_{ev, *}^{w}$ of all weak Jacobi forms of even weight is a
polynomial algebra in two variables over $\frak
M_*(SL_2(\bz))=\oplus_k \frak M_{k}(SL_2(\bz))$ with the standard
generators $\phi_{0,1}$ and $\phi_{-2,1}$:
$$
\align
\phi_{-2,1}(\tau, z)&
=\biggl(\frac{\vth(\tau,z)}{\eta(\tau)^3}\biggl)^2
=\biggl(r^{-1/2}\prod_{n\ge 1}\,(1-q^{n-1} r) (1-q^n
r^{-1})(1-q^n)^{-2}\biggr)^2
\\
{}&=(r-2+r^{-1})-2q(r^2-4r+6-4r^{-1}+r^{-2})+q^2(\dots),
\tag{3.5}
\\
\phi_{0,1}(\tau, z)&=\Delta(\tau)^{-1}\phi_{12,1}(\tau,z)\\
{}&=(r+10+r^{-1})+q(10r^{-2}-64r^{-1}+108-64r+10r^2)+q^2(\dots ),
\endalign
$$
where $\vth$ is the Jacobi theta-series \thetag{1.4} and $\phi_{12,1}$ is
the unique Jacobi cusp form of weight $12$ and index $1$ with integral
coprime Fourier coefficients and $q=\exp{(2\pi i \tau)}$,
$r=\exp{(2\pi i z)}$.

We can represent $\Delta^{-N}\phi_{2k,m}$ as a polynomial in
$\phi_{-2,1}$ and $\phi_{0,1}$ over $\frak M_{*}(SL_2(\bz))$:
$$
\Delta^{-N}\phi_{2k,m}=
\phi_{-2,1}^{6N-k}\phi_{0, m-6N+k}=
\phi_{-2,1}^{6N-k}\sum\Sb
0\le i\le m-6N+k\\
\vspace{0.5\jot}
i\ne 1
\endSb
f_{2i}\phi_{-2,1}^i \phi_{0,1}^{m-6N+k-i}
\tag{3.6}
$$
where the $f_{2i}$ are modular forms of weight $2i$ with respect to
$SL_2(\bz)$.  This shows $6N\ge k+m.$ If $6N = k+m$, then the Jacobi
form
$$
\phi_{2k,m}=\Delta^N\phi_{-2,1}^m =\eta^{6(4N-m)}\vth^{2m}
$$
is of the type stated in the theorem.

Let $6N> k+m$.  The Fourier expansion of the weak Jacobi form
$\phi_{0, m-6N+k}$ has a non-zero $q^0$-term due to the choice of
$N=\hbox{ord}_q(\phi_{2k,m})$.  Let us for the moment suppose that the
$q^0$-term of $\phi_{0, m-6N+k}$ is not a constant.  Due to the
relation $\phi_{2k,m}(\tau, -z)=\phi_{2k,m}(\tau, z)$, the $q^0$-term
of the weak Jacobi form $\phi_{0, m-6N+k}$ contains a summand of
type $a(r^{\beta}+r^{-\beta})$ with $a\ne 0$ and a positive integer
$\beta$ which we assume to be maximal.  This implies that the Fourier
expansion of the Jacobi cusp $\phi_{2k,
m}=\Delta^N(\phi_{-2,1})^{6N-k}\phi_{0, m-6N+k}$ has a non-zero
Fourier coefficient of type $q^Nr^{6N-k+\beta}$.  It follows that the
norm of the index of this coefficient has to be positive, i.e. $
4Nm-(6N-k+\beta)^2>0 $.  Therefore, since $\beta\ge 1$, we have
$$
N< \frac{1}{18}\biggl(
 (3k-3\beta+m)+\sqrt{(3k-3\beta+m)^2-9(k-\beta)^2}\biggr)\le
\frac{3k-3+m}9.
$$
Thus the equality in \thetag{3.2} can hold only if the minimum equals
$\frac {k+m}6$. This case was considered above.  To finish the proof
of Proposition 3.2 it remains to prove the assumption about the
$q^0$-part of $\phi_{0,t}$ if $t\ne 0$.
\enddemo
\proclaim{Lemma 3.4}Let $\phi_{0,t}$ be a weak Jacobi form of weight
zero and index $t$ (\,$t\ne 0$) with $\hbox{\rm ord}_q(\phi_{0,t})=0$.
Then the $q^0$-part of the Fourier expansion of $\phi_{0,t}$ is not a
constant.
\endproclaim
\demo{Proof of the lemma}
Similar to \thetag{3.6} we obtain
$$
\phi_{0,t}(\tau,z)=\sum\Sb
0\le i\le t\\
\vspace{0.5\jot}
i\ne 1
\endSb
f_{2i}(\tau)\phi_{-2,1}^i(\tau, z) \phi_{0,1}^{t-i}(\tau, z)
\tag{3.7}
$$
where $f_{2i}(\tau)=a_iE_{2i}(\tau)+\dots$ ($a_i\in \bc$) are modular
forms from $\frak M_{2i}(SL_2(\bz))$ and $E_{2i}(\tau)$ are the
Eisenstein series of weight $2i$ with respect to $SL_2(\bz)$ with the
constant term $1$. We can assume that the other summands are cusp forms.

According to \thetag{3.5} the $q^0$-term $\phi_{0,t}^{(0)}(z)$ of
$\phi_{0,t}(\tau,z)$ is equal to
$$
\phi_{0,t}^{(0)}(z)=\sum\Sb
0\le i\le t\\
\vspace{0.5\jot}
i\ne 1
\endSb
a_i\, X^i (X+12)^{t-i} \qquad\quad(X=r-2+r^{-1}).
$$
If $\phi_{0,t}^{(0)}$ is a constant, then we obtain a triangular
linear system of $t$ equations with $t$ unknowns $a_i$ ($i\ne 1$!)
$$
X^k\sum_{i\le k} a_i C_{t-i}^{t-k}=0 \qquad\quad(1\le k\le t).
$$
This implies that all $a_i=0$ and hence $\phi_{0,t}^{(0)}= 0$, a
contradiction to the assumption about the $q$-order of $\phi_{0,t}$.

\hfill\hfill\qed
\enddemo
\demo{Proof of Theorem 3.1}
The ``only if\," part of the first claim of the theorem follows from
Corollary 1.6.

Let $t=1,\ 2,\ 4,\ 5$.  According to \thetag{1.15} we have to prove
that there are no cusp forms of weight $3$ with respect to
$\Gamma_t'$.  For this it is clearly enough to prove that for any
character $\chi_{a,b}:\Gamma_t\to \bc^*$ (see Notation 2.2)  of
the paramodular group $\Gamma_t$ for $t=2,\ 4, \ 5$ we have
$\hbox{dim}\, (\frak N_3(\Gamma_t, \chi_{a,b}))=0$.  (For $t=1$ this
follows from Igusa's result \cite{Ig2} about the graded ring of Siegel
modular forms for $Sp_4(\bz)$, but it can be easily obtained by the
same method which we apply for $t=5$.)  For $t=3,\ 6,\,7$ we will
show that the corresponding spaces of cusp forms of weight $3$ contain
only one
function.

Let $F_{\chi_{a,b}}\in \frak N_k(\Gamma_t, \chi_{a,b})$.  The
restriction of $\chi_{a,b}$ to the centre of the maximal parabolic
subgroup $\Gamma_{\infty,t}$, which is generated by the element $B$
(see \thetag{2.1}), determines the type of the Fourier-Jacobi
expansion of the cusp form $F$.  Let $t_2=(2t, 12)$ as in Theorem 2.1.
Then
$$
F_{\chi_{a,b}}(Z)=
\sum\Sb m\equiv b\, mod\, t_2\\
\vspace{0.5\jot}
m>0\endSb f_{k, tm/t_{2}}^{(a)}(\tau, z)\exp(2\pi i m\frac{t}{t_2}
\omega).
\tag{3.8a}
$$
The character of the Jacobi forms $f_{k, tm/t_{2}}$ is equal to
$v_\eta^{24a/t_2}\times v_H^{\varepsilon}$, where $\varepsilon \equiv
2tb/t_2\ \hbox{mod}\ 2$.  We also have $ F|_k V_t\in \frak
N_k(\Gamma_t, \chi_{b,a}) $ since $V_tBV_t=A$ and the character
$\chi_{b,a}$ is $V_t$-conjugate to $\chi_{a,b}$, i.e.
$\chi_{b,a}(\gamma)=\chi_{a,b}(V_t \gamma V_t)$ for $\gamma\in
\Gamma_t$.  Therefore, we can define a $V_t$-invariant modular form
$F\cdot F|_k V_t\in \frak N_{2k}(\Gamma_t, \chi_{a+b,a+b})$.  It has
the following Fourier-Jacobi expansion
$$
(F\cdot F|_k V_t)(Z)=
\sum\Sb m\equiv\, a+b\, mod\, t_2\\
\vspace{0.5\jot}
m\ge a+b \endSb f_{2k, tm/t_{2}}(\tau, z)\exp(2\pi i m\frac{t}{t_2}
\omega).
\tag{3.8b}
$$
We divide the proof of Theorem 3.1 into six steps depending on the
values of the polarization $t$.
\smallskip

({\bf V}). Let $t=5$.  According to Theorem 2.1 there is only one
non-trivial character $\chi_2=\chi_{1,1}:\Gamma_5\to \{\pm 1\}$.
(If the character is trivial our result follows from Corollary 3.3.)  Let
$F\in \frak N_3(\Gamma_5, \chi_{1,1})$.  Then $G^{(5)}=F\cdot F|V_5\in
\frak N_6^+(\Gamma_5)$ is a $V_5$-invariant cusp form.
According to
the dimension formula \thetag{3.4} $\hbox{dim}( J_{6,5}^c)=1$ and a
basis of this space can be given explicitly in terms of Jacobi forms
of half-integral indices, namely
$$
J_{6,5}^c=\bc\phi_{6,5}, \qquad\text {where}\quad
\phi_{6,5}=\eta^2\vth^6[\vth, \vth_{3/2}]
$$
(see \thetag{1.7} and \thetag{1.12}).  Calculating the Fourier
expansion of $[\vth, \vth_{3/2}]$ one gets
$$
\frac{1}2\eta^{-4}[\vth, \vth_{3/2}]=(r+4+r^{-1})+q(\dots).
$$
Thus $\phi_{6,5}$ {\it is not a square of a Jacobi form of weight $3$
and index $5/2$}.
It follows that the first Fourier-Jacobi coefficient $g_{6, 5M}$ of
$G^{(5)}$ has $M\ge 2$.  The $q$-order of $g_{6,5M}$ is at least $M$,
since the $V_5$-invariant cusp form $G^{(5)}$ has order $5M$ with
respect to the variable $\omega$. We shall use this type of argument
repeatedly in the sequel.  According to Proposition 3.2
$\hbox{ord}_q(g_{6,5m})\le \frac{6+5m}9<m$ if $m\ge 2$.  Hence
$G^{(5)}=0$.
\smallskip

({\bf II}) and ({\bf IV}).  Let $t=2$ or $t=4$.  An arbitrary
character $\chi_{a,b}$ of $\Gamma_t$ has order $1$, $2$ or $4$. Let us
consider cusp forms $F=F_{\chi_{a,b}}\in \frak N_3(\Gamma_t,
\chi_{a,b})$. Then
$$
G^{(t)}(Z)=\bigl((F\cdot F|V_t)(Z)\bigr)^4=\sum_{m\ge 2} g_{24, tm}^{(t)}(\tau ,
z)\exp(2\pi i tm\omega)\in \frak N_{24}^+(\Gamma_t)
\qquad(t=2,\,4)
$$
and again use Proposition 3.2. The inequality
$\hbox{ord}_q(g_{24,tm})\le\frac{12+tm}6<m$ holds if $m\ge 4$ for
$t=2$ or $m\ge 7$ for $t=4$.  Thus to prove that $G^{(2)}=0$,
resp. $G^{(4)}=0$, it is enough to show that the Fourier-Jacobi
expansion of $F\cdot F|V_t$ starts with coefficients $f_{6,2}$ or $f_{6,7}$ for
$t=2$ or $4$ respectively.  To show this we check that for any $b\le a$
such that $a+b<4$, resp. $a+b<7$ (see \thetag{3.8b}) the first
possible Fourier-Jacobi coefficient $f_{k,bt/4}^{(a)}$ of any
$F_{\chi_{a,b}}\in \frak N_3(\Gamma_t, \chi_{a,b})$ is zero.  For this
we have to consider the character $\chi_{1,1}$ of the group $\Gamma_2$
and the characters $\chi_{a,1}$ ($1\le a \le 4$), $\chi_{a,2}$ ($2\le
a \le 4$) and $\chi_{3,3}$ of $\Gamma_4$.  We have
$$
f_{3,1/2}^{(1)}\cdot \eta^3\vth^5\in J_{7,3}=\{0\},
\qquad
f_{3, b}^{(a)}\cdot \vth^{8-2a}\in J_{7-a,b+4-a}^c=\{0\}
$$
for the values of $a$ and $b$ chosen above (see \thetag{3.4}).
\smallskip

({\bf III}). Let $t=3$. By Theorem 2.1 the order of any character of
$\Gamma_3$ is a divisor of $6$.  For given $F=F_{\chi_{a,b}}\in \frak
N_3(\Gamma_3, \chi_{a,b})$ we consider the form
$$
G^{(3)}(Z)=(F\cdot F|V_3)(Z)^6=\sum_{m\ge 2} g_{36, 3m}(\tau ,
z)\exp(6\pi i m\omega)\in \frak N_{36}^{+}(\Gamma_3).
$$
According to Proposition 3.2 $\hbox{ord}_q(g_{36,3m})\le\frac{18+3m}6<
m$ if $m\ge 7$.  Thus if we prove that the Fourier-Jacobi coefficients
of $F\cdot F|V_3$ of indices smaller than $7/2$ are zero, then it
follows that $F_{\chi_{a,b}}=0$.

Recall that $\chi_{a,b}(A)^3=\chi_{a,b}(B)^3$ (see Theorem 2.1).
Therefore, due to \thetag{3.8a} and $\thetag{3.8b}$ we have to consider
only the characters $\chi_{1,1}$, $\chi_{3,1}$, $\chi_{5,1}$,
$\chi_{2,2}$, $\chi_{4,2}$ and $\chi_{3,3}$.  For the second, the
third, the fourth and the fifth character we have that the first
Fourier-Jacobi coefficient in the expansion
\thetag{3.8a} of the  corresponding cusp form $F_{\chi_{a,b}}$
is zero:
$$
\align
f_{3,1/2}^{(3)}\cdot \eta^3\vth^3&\in J_{6,2}^c=\{0\} \quad
(\chi_{3,1}),
\qquad
f_{3,1/2}^{(5)}\cdot \eta\vth\in J_{4,1}^c=\{0\} \quad (\chi_{5,1}),\\
f_{3,1}^{(2)}\cdot \eta^4\vth^4&\in J_{7,3}=\{0\} \quad (\chi_{2,2}),
\qquad
f_{3,1}^{(4)}\cdot \eta^2\vth^2\in J_{5,2}=\{0\} \quad (\chi_{4,2}).
\endalign
$$
Thus the cusp forms $F_{\chi_{a,b}}$ of weight $3$ are zero for these
characters.

For the first character the situation is more complicated.  The first
Fourier-Jacobi coefficient of $F_{\chi_{1,1}}$ satisfies the relation
$$
f_{3,1/2}^{(1)}\cdot\eta^5\vth^5\in J_{8,3}^c=
\bc\cdot  \eta^6\vth^4 \phi_{3,1},
$$
where $\phi_{3,1}(\tau,z)=\phi_{12,1}(\tau,z)/\eta(\tau)^{18}
\in J_{3,1}(v_\eta^{6}\times 1_H)$
and $\phi_{12,1}$ is the unique (up to a constant) Jacobi cusp form of
weight $12$ and index $1$.  The Jacobi form $\phi_{3,1}(\tau,z)$ does
not vanish for $z=0$ since $\phi_{12,1}(\tau,0)=\eta(\tau)^{24}$.
Thus $f_{3,1/2}^{(1)}=0$.  As a result we have proved that there are
no cusp forms of weight $3$ for all but 1 character of the group
$\Gamma_3$:
$$
\bigoplus_{(a,b)\ne (3,3)} \frak N_3(\Gamma_3, \chi_{a,b})=\{0\}.
\tag{3.9}
$$

Let us consider the last character $\chi_{3,3}=\chi_2$
(see \thetag{1.10a-b}) of order $2$.
This character can
be extended to a character of the group $\Gamma_3^+$, thus we can
consider a decomposition of the space
$\frak N_3(\Gamma_3, \chi_{2})$ of type \thetag{3.3},
namely $\frak N_3(\Gamma_3, \chi_{2})=
\frak N_3(\Gamma_3^+, \chi_{3,3,+})\oplus
\frak N_3(\Gamma_3^+, \chi_{3,3,-})$.
The first Fourier-Jacobi coefficient of any  element
$F_{\chi_{2}}\in \frak N_3(\Gamma_3^+, \chi_{3,3,\pm})$
is a Jacobi form of index $1/2$
$$
f_{3,1/2}^{(3)}\in J_{3,\frac{3}2}^c(v_{\eta}^{12}\times v_H)=
\bc \cdot \eta^3\vth^3.
$$
The last space is one-dimensional since
$\bigl(J_{3,\frac{3}2}^c(v_{\eta}^{12}\times v_H)\bigr)^2
\subset J_{6,3}^c$ and $\hbox{dim\,}J_{6,3}^c=1$.
We can define the lifting of $\eta\vth$ (see \thetag{1.9})
$$
\Delta_1(Z)=\hbox{Lift}(\eta\vth)(Z)\in \frak N_1(\Gamma_3, \chi_{1,1}).
$$
Thus there exists a constant $c$ such that
$$
F_{\chi_{2}}(Z)^2- c \Delta_1(Z)^{6}=
\sum_{m\ge 2} f_{6, 3m}(\tau, z) e^{6\pi i m\omega}\equiv 0
$$
according to Proposition 3.2.
\smallskip

({\bf VII}). Let $t=7$. According to Theorem 2.1 there exists only one
non-trivial character $\chi_2=\chi_{1,1}$ of $\Gamma_7$ of order $2$.
(The case of trivial character was considered in Corollary 3.3.)
We remark that
$$
J_{3,\frac{7}2}^c(v_{\eta}^{12}\times v_H)=
\bc\cdot \vth^2[\vth, \vth_2]
$$
(see \thetag{1.12}). This space is spanned by one Jacobi form, since
for any $\phi\in J_{3,\frac{7}2}^c(v_{\eta}^{12}\times v_H)$ we obtain
$\phi(\eta\vth)^3\in J_{6, 5}^c$, and the last space is one
dimensional.

Let
$$
F^{(7)}(Z)=
\sum \Sb m\equiv 1  \,mod\, 2\endSb
f_{3, \frac{7m}2}(\tau,z)e^{7\pi i m \omega}
\in \frak N_3(\Gamma_7^+, \chi_{2,\pm}).
$$
Then there exists a constant $c$ such that
$$
F^{(7)}(Z)^2- c\, \Xi(Z)^2=
\sum_{m\ge 2} f_{6, 7m}(\tau,z) e^{14\pi i m\omega}
$$
where $\Xi=\hbox{Lift}(\vth^2[\vth, \vth_2])$ (see \thetag{t=7} in the
proof of Theorem 1.4) and $f_{6, 7m}\in J_{6, 7m}^c$.  Due to
Proposition 3.2 $
\hbox{ord}_q(f_{6,7m})<\frac{6+7m}9\le m
$ if $m\ge 3$.  Thus $(F^{(7)})^2=c \,\Xi^2$ if we prove that there is no
Jacobi cusp form $f_{6,14}$ of weight $6$ and index $14$ such that
$\hbox{ord}_q(f_{6,14})\ge 2$.  We shall now show this.

If $\Delta^{-3}f_{6,14}$ is a weak holomorphic Jacobi form, then using
representation \thetag{3.7} we obtain
$\Delta^{-3}f_{6,14}=(\phi_{-2,1})^{15}\phi_{0,-1}$, a contradiction.
Thus $f_{6,14}$ would have $q$-order 2 and
$$
\Delta^{-2}f_{6,14}
=(\phi_{-2,1})^9\phi_{0,5}, \qquad \phi_{0,5}\in J_{0,5}^{weak}.
$$
If $\phi_{0,5}$ contains $r^d$ with $d\ge 2$, then we get a
contradiction as in Lemma 3.4.  For $t=5$ all modular forms in the
representation
\thetag{3.7} are Eisenstein series.
Thus, there is only one weak Jacobi form $\phi_{0,5}(\tau,z)$ with
$q^0$-term equal to $r^{1}+c+r^{-1}$.  It is easy to determine this
function using $\thetag{3.7}$ or using two weak Jacobi forms found in
\cite{GN, Example 4.5}.  As a result we obtain
$$
5\phi_{0,5}(\tau,z)=(5r^1+2+5r^{-1})+q(-r^5+\dots).
$$
The Jacobi form $\eta^{48}(\phi_{-2,1})^{9}\phi_{0,5}$ is not
holomorphic, since its Fourier expansion contains the term
$cq^3r^{14}$ $(c\ne 0)$. Therefore $f_{6,14}=0$.  This finishes the
proof for $t=7$.
\smallskip

({\bf VI}). Let $t=6$.  Any cusp form $F_{\chi_2}$ of weight $3$ with
character $\chi_2=\chi_{6,6}$ of order $2$ has a Fourier expansion
$$
F_{\chi_2}(Z)=\sum
\Sb m\equiv 1 \,mod\,2\endSb
f_{3, {3m}}(\tau,z)e^{6\pi i m\omega}\in \frak N_3^{\pm}(\Gamma_6,
\chi_{2}).
$$
Similar to (III) (see also \thetag{t=6} of the proof of Theorem 1.1)
$$
J_{3,3}^{c}(v_\eta^{12}\times 1_H)=
\bc\cdot \eta^3\vth^2\vth_2
$$
is one-dimensional, since $\hbox{dim}\,J_{6,6}^c=1$.  Thus there
exists a constant $c$ such that
$$
F_{\chi_2}(Z)^2- c\al(\eta^3\vth^2\vth_2)(Z)^2=
\sum_{m\ge 2} f_{6, 6m}(\tau,z) e^{12\pi i m\omega}\equiv 0,
$$
since due to Proposition 3.2 $
\hbox{ord}_q(f_{6,6m})< \frac{6+6m}9\le m
$ for $m\ge 2$.
Hence $h^{3,0}(\widetilde{\Cal A}_6(\chi_{2}))=1$.
\medskip

We remark that we have in fact proven that
$\hbox{dim}(\frak N_6^+(\Gamma_6))=1$.
The same arguments show us that
$\hbox{dim}(\frak N_6^-(\Gamma_6))=0$ since the generator
$(\eta^3\vth^2\vth_2)^2$ of the space $J_{6,6}^c$ has
$q$-order one.

Now we can find the geometric genus of
$\widetilde{\Cal A}_6(\chi_{3})$ where
$\chi_3=\chi_{4,4}:\Gamma_6\to \bc^*$ is the character induced by
$v_\eta^8\times 1_H$ (see \thetag{1.10a--b} and Notation 2.2).
Let $\chi_{a,b}$ be  a different character of $\Gamma_6$.
It follows from the definition of the characters $\chi_{a,b}$ that
$\hbox{ker}(\chi_{4,4})=\hbox{ker}(\chi_{a,b})$ if and only if
$\chi_{a,b}=\chi_{4,4}^2=\chi_{8,8}$.
In Theorem 1.4 we constructed the cusp form
$\al(\eta^5\vth_2)\in \frak N_3(\Gamma_6, \chi_{4,4})$.
A modular  form from $\frak N_3(\Gamma_6, \chi_{4,4})$ was found in
\cite{GN, Example 1.15}. This is the $2$-lifting (or $-1$-lifting)
of the Jacobi form $\eta^5\vth_2$:
$$
\al_2(\eta^5\vth_2)=\sum\Sb m\equiv 2\,mod\,3\endSb
m^{-1}(\eta(\tau)^5\vth(\tau, 2z)\exp(4\pi i \omega))|_3 T_-^{(3)}(m)
\in \frak N_3(\Gamma_6, \chi_{8,8}).
$$
Let $F_{\chi_{4,4}}\in \frak N_3(\Gamma_6, \chi_{4,4})$ and
$F_{\chi_{8,8}}\in \frak N_3(\Gamma_6, \chi_{8,8})$ be arbitrary cusp forms.
Then it is clear that
$F_{\chi_{4,4}}\cdot F_{\chi_{8,8}}\in \frak N_6(\Gamma_6)$.
This space is one dimensional (see above), thus
$$
h^{3,0}(\widetilde{\Cal A}_6(\chi_{4,4}))=
\hbox{dim}(\frak N_3(\Gamma_6, \chi_{4,4}))+
\hbox{dim}(\frak N_3(\Gamma_6, \chi_{8,8}))=2.
$$
Hence Theorem 3.1 is proved.
\hfill\hfill\qed
\enddemo
\medskip
We want to conclude this section with a brief discussion of the commutative
neighbours of ${\Cal A}_3$. This case is of special geometric interest and is
also closely related to the Calabi-Yau threefold found by Barth and Nieto
\cite{BN}.

\proclaim{Corollary 3.5}
There are exactly 12 Siegel modular varieties between ${\Cal A}_3^{com}$ and
${\Cal A}_3$. Of these six have geometric genus $1$,
the others have geometric genus $0$.
\endproclaim
\demo{Proof} Recall from Theorem 2.1 that
$\Gamma_3/\Gamma'_3\cong\bz/3\times\bz/6$. This group is generated by elements
$A$ and $B$ with the relations $A^6=B^6=1$ and $A^3=B^3$. The Siegel modular
varieties between ${\Cal A}_3^{com}$ and ${\Cal A}_3$ are in 1:1 correspondence
with the subgroups of $G_3=\Gamma_3/\Gamma'_3$.
The lattice of subgroups of $G_3$ is represented in Diagram 1.
The proof of Theorem 3.1 shows that for a group $\Gamma$ between $\Gamma'_3$
and $\Gamma_3$ the corresponding moduli space ${\Cal A}(\Gamma)$ has geometric
genus 1 if and only if $\Gamma$  is contained in
$\operatorname{ker}(\chi_{3,3})$. The
latter is generated by $A^2$ and $B^2$ and isomorphic to $\bz/3\times \bz/3$.
This, togehter with the diagram below shows the claim.

\hfill\hfill{\qed}
\enddemo

\newpage
\midinsert
$$
\vbox{\centerline{\epsfxsize=5in\epsfbox{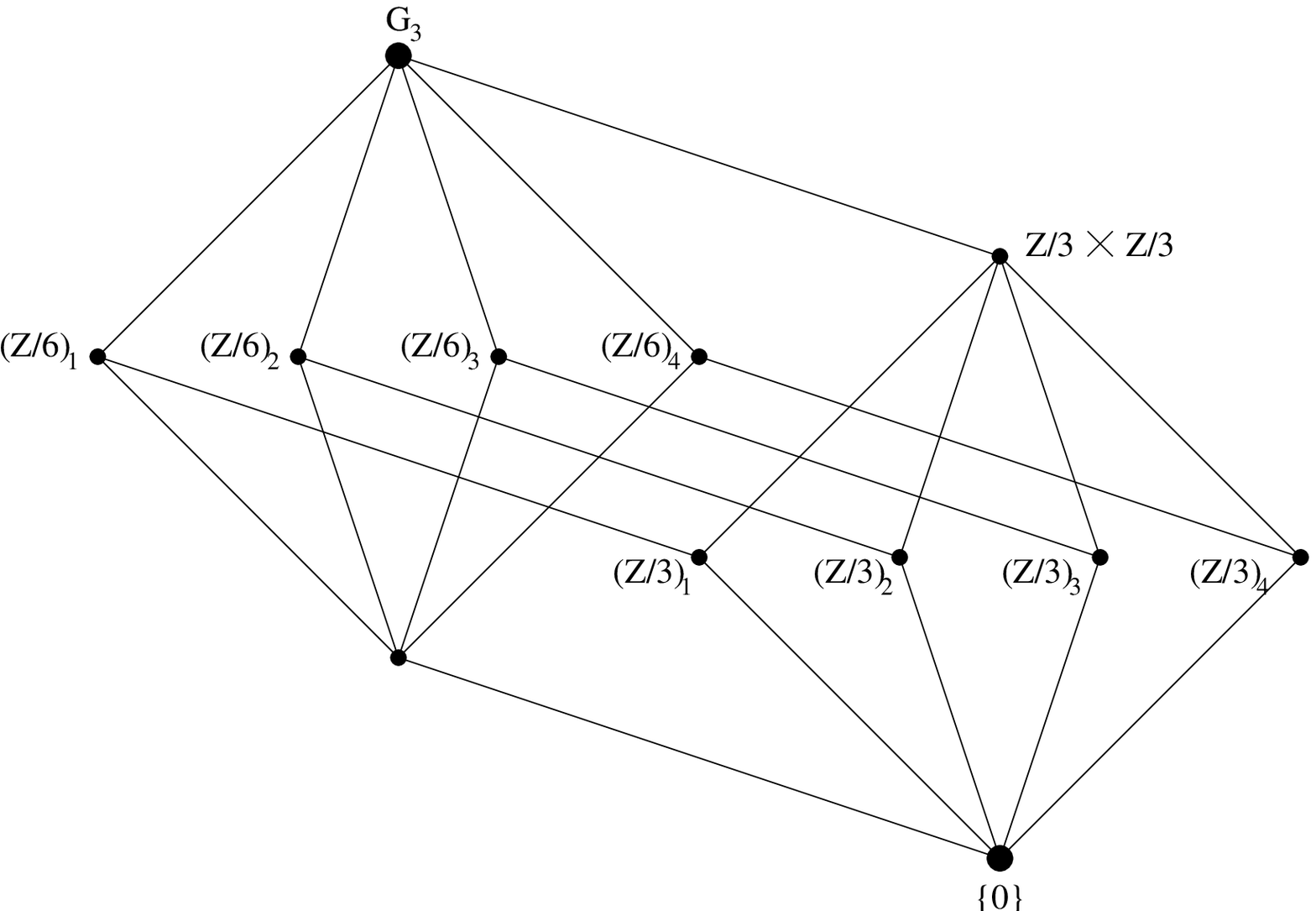}}
\vskip20pt
\centerline{
{\bf Diagram 1.} Lattice of the subgroups of the group
$G_3=\Gamma_3/\Gamma_3'$.}}
$$
\endinsert

\remark{Remark 3.6}
The lattice of the subgroups of $G_3$ corresponds to a diagram of coverings
of Siegel modular varieties. All but one of these coverings are of type
$\Cal A(\chi)$ for some character $\chi$ of the paramodular group $\Gamma_t$.
To be more precise we have to label the subgroups
$(\bz/3)_i$ and $(\bz/6)_i$. We put
$$
\aligned
(\bz/3)_1&=<A^2>, \qquad\\
(\bz/3)_3&=<A^2B^2>, \qquad
\endaligned
\aligned
(\bz/3)_2&=<B^2>,\\
(\bz/3)_3&=<A^4B^2>=<AB^{-1}>=<A^{-1}B>.
\endaligned
$$
This also determines the subgroups $(\bz/6)_i$.
Note that
$$
\bz/3\times \bz/3=<A,B>, \qquad\bz/2=<A^3>=<B^3>.
$$
Altogether we have $18$ characters $\chi_{a,b}$;
$1\le a\le b$, $a-b\equiv 0 \pmod 2$. Each of the subgroups $(\bz/3)_i$,
resp. $(\bz/6)_i$ is the kernel of two different characters, the subgroup
$\bz/3\times \bz/3$ is the kernel of $\chi_{3,3}$.
The precise relation between subgroups of $G_3$ and characters is given
by the following table
$$
\spreadmatrixlines{1\jot}
\matrix
\format
\c&\qquad\quad \c &\qquad \c\\
[\,\hbox{{group}}\,] & [\,\hbox{{character}}\,]
&[\,\hbox{{ order of character}}\,]\\
\vspace{2\jot}
(\bz/3)_1 &\chi_{3,1}\ \text{ and }\ \chi_{3,5} &6\\
(\bz/3)_2 &\chi_{1,3}\ \text{ and }\ \chi_{5,3} &6\\
(\bz/3)_3 &\chi_{1,5}\ \text{ and }\ \chi_{5,1} &6\\
(\bz/3)_4 &\chi_{1,1}\ \text{ and }\ \chi_{5,5} &6\\
(\bz/6)_1 &\chi_{6,2}\ \text{ and }\ \chi_{6,4} &3\\
(\bz/6)_2 &\chi_{2,6}\ \text{ and }\ \chi_{4,6} &3\\
(\bz/6)_3 &\chi_{4,2}\ \text{ and }\ \chi_{2,4} &3\\
(\bz/6)_4 &\chi_{2,2}\ \text{ and }\ \chi_{4,4} &3\\
\bz/3\times \bz/3 &\chi_{3,3}                   &2\\
G_3 &\chi_{6,6}&1
\endmatrix
$$
The subgroup $\bz/2$ is not the kernel of a character, but it can be
written as the intersection of the kernel of two characters in several ways.
Note that for example
$(\bz/6)_1=<A>$, $\ (\bz/6)_2=<B>$ and hence
$\bz/2=\hbox{ker\,}(\chi_{6,2})\cap \hbox{ker\,}(\chi_{2,6})$.

The lattice of the subgroups of $G_3$ then corresponds to the
following diagram of coverings

\midinsert
$$
\vbox{\centerline{\epsfxsize=5in\epsfbox{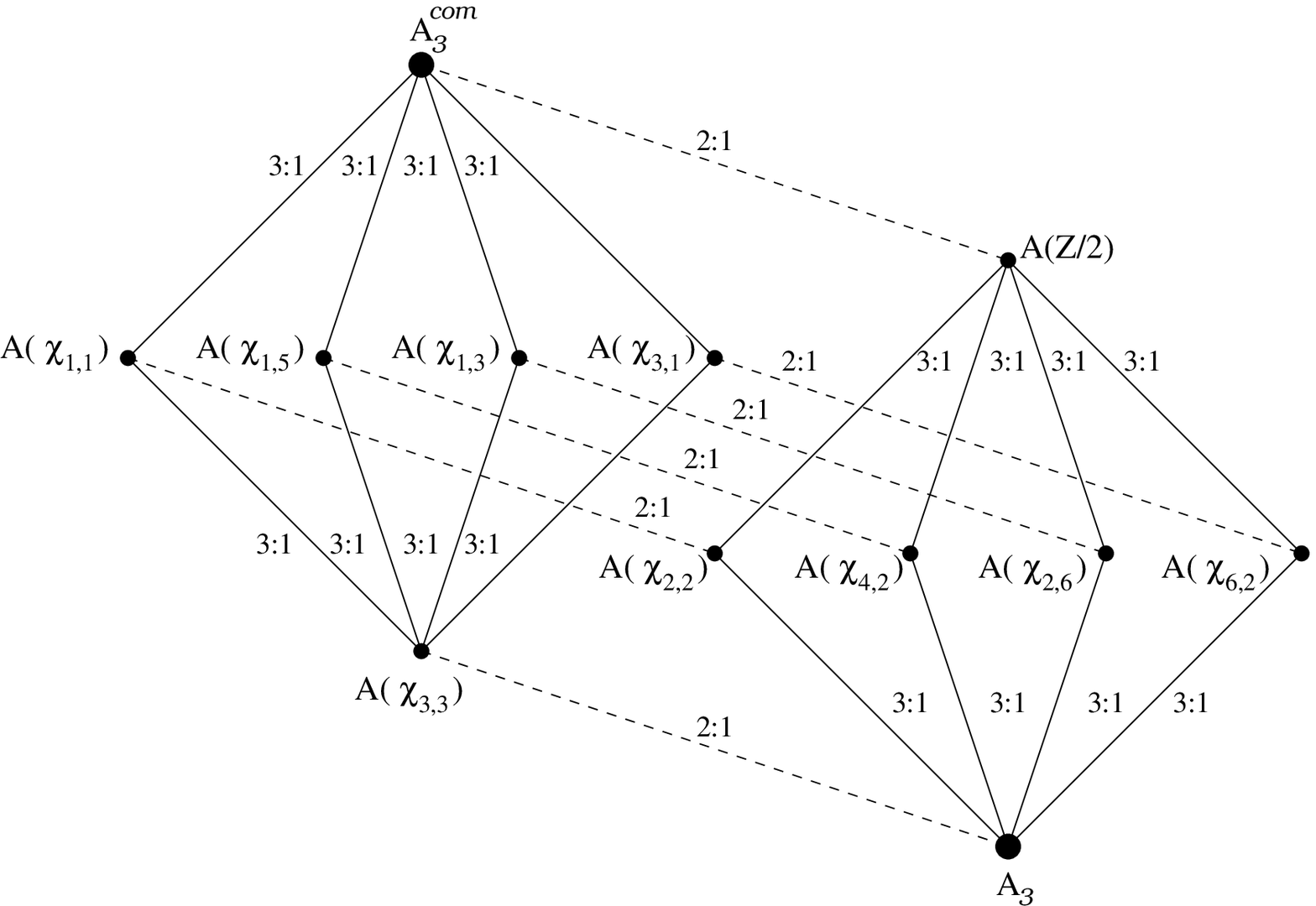}}
\vskip20pt
\centerline{
{\bf Diagram 2.} Commutative coverings of $\Cal A_3$.}}
$$
\endinsert

The modular varieties corresponding to the vertices of
the left square, resp.  right square have geometrical genus $1$,
resp. $0$.

\endremark

\remark{Remark 3.7}
There are exactly 20 Siegel modular varieties between ${\Cal A}_3^{com}$ and
${\Cal A}^+_3=\Gamma^+_3\backslash \bh_2$.
(Not all of these are abelian covers of ${\Cal A}^+_3$.)
Of these, exactly $6$ have geometric genus $1$, whereas the
others have geometric genus $0$ (by \thetag{3.9}). 
The six modular varieties of geometric
genus $1$ are natural candidates for Calabi--Yau varieties, whereas we
expect the other varieties to be unirational.

\endremark

%\newpage

\enddocument
\end